\documentclass[a4paper,12pt]{article}     

\usepackage{graphicx}		  
\usepackage{amsmath,amssymb}	

\newtheorem{lemma}{Lemma}[section]
\newtheorem{theorem}[lemma]{Theorem}

\newtheorem{assumption}{Assumption}

\def\authorfont{\footnotesize}

\def\ccode#1{\par
\vspace*{8pt}
{\authorfont{\leftskip18pt\rightskip\leftskip
\noindent #1\par}}\par}

\newenvironment{Proof}{
\hspace*{-9mm}
{ \it Proof.}}
{\hfill {$\square$}\vspace{1.5em}}

\begin{document}

\begin{center}{
{\Large 
 Properties of minimal charts and
 their applications V: 
 charts of type $(3,2,2)$}
\vspace{10pt}
\\ 
Teruo NAGASE and Akiko SHIMA
}
\end{center}

\date{2019.1.22}

\begin{abstract}
Let $\Gamma$ be a chart,
and we denote by $\Gamma_m$
the union of all the edges of label $m$.
A chart $\Gamma$ is of type $(3,2,2)$
if there exists a label $m$
such that 
$w(\Gamma)=7$,
$w(\Gamma_m\cap\Gamma_{m+1})=3$,
$w(\Gamma_{m+1}\cap\Gamma_{m+2})=2$,
and $w(\Gamma_{m+2}\cap\Gamma_{m+3})=2$
where 
$w(G)$ is the number of white vertices in $G$.
In this paper, we prove that there is 
no minimal chart of 
type $(3,2,2)$.
\end{abstract}

%
%
%
%

\ccode{2010 Mathematics Subject Classification. Primary 57Q45; Secondary 57Q35.}
\ccode{ {\it Key Words and Phrases}. surface link, chart, white vertex. }

\setcounter{section}{0}
\section{Introduction}


Charts are oriented labeled graphs in a disk (see  \cite{KnottedSurfaces},\cite{BraidBook}, and see Section~\ref{s:Prel}  for the precise definition of charts).
From a chart, we can construct an oriented closed surface 
embedded in 4-space ${\Bbb R}^4$ 
 (see \cite[Chapter 14, Chapter 18 and Chapter 23]{BraidBook}). 
A C-move 
is a local modification between two charts
in a disk (see Section~\ref{s:Prel} for C-moves).
A C-move between two charts induces 
an ambient isotopy between oriented closed surfaces 
corresponding to the two charts.

We will work in the PL category or smooth category. All submanifolds are assumed to be locally flat.
In \cite{ONS},
we showed that there is no minimal chart with exactly five vertices
 (see Section~\ref{s:Prel} for the precise definition of minimal charts). 
Hasegawa proved that there exists a minimal chart with exactly
six white vertices \cite{H1}. 
This chart represents a 2-twist spun trefoil.
In \cite{INS} and \cite{NST},
we investigated minimal charts with exactly four white vertices.
In this paper, 
we investigate properties of minimal charts and
need to prove that
there is no minimal chart with exactly seven white vertices
(see \cite{ChartApp1},\cite{ChartAppII},
\cite{ChartAppIII},\cite{ChartAppIV},
\cite{ChartAppVI}).

Let $\Gamma$ be a chart.
For each label $m$, we denote by $\Gamma_m$
the union of all the edges of label $m$.

Now we define a type of a chart:
Let $\Gamma$ be a chart, 
and $n_1,n_2,\dots,n_p$ integers.
The chart $\Gamma$ is of {\it type $(n_1,n_2,\dots,n_k)$} if there exists a label $m$ of $\Gamma$ satisfying the following three conditions:
\begin{enumerate}
\item[(i)] For each $i=1,2,\dots, k$, 
the chart $\Gamma$ contains exactly $n_{i}$ white vertices in $\Gamma_{m+i-1}\cap \Gamma_{m+i}$.
\item[(ii)] If $i<0$ or $i>k$, then $\Gamma_{m+i}$ does not contain any white vertices.
\item[(iii)] Both of the two subgraphs $\Gamma_m$ and $\Gamma_{m+k}$ contain at least one white vertex.
\end{enumerate}
If we want to emphasize the label $m$,
then we say that $\Gamma$ is of {\it type $(m;n_1,n_2,\dots,n_k)$}. 
Note that $n_1\ge1$ and $n_k\ge1$ by the condition (iii).

We proved in \cite[Theorem 1.1]{ChartAppII} that
if there exists a minimal $n$-chart $\Gamma$ with exactly seven white vertices,
then $\Gamma$ is a chart of 
type $(7),(5,2),(4,3),(3,2,2)$ or $(2,3,2)$ 
(if necessary we change the label
$i$ by $n-i$ for all label $i$).
The following is the main result in this paper.

\begin{theorem}
\label{MainTheorem} 
There is 
no minimal chart of 
type $(3,2,2)$.
\end{theorem}

The paper is organized as follows.
In Section~\ref{s:Prel},
we define charts and minimal charts.
In Section~\ref{s:ConnectedComponent},
we investigate connected components of $\Gamma_m$
with at most three white vertices
for a minimal chart $\Gamma$.
In Section~\ref{s:KAngledDisk},
we review a $k$-angled disk,
a disk whose boundary consists of edges of label $m$ and contains exactly $k$ white vertices.
In Section~\ref{s:DiskTwo},
we investigate
a disk $D$ with exactly two white vertices 
of $\Gamma_m$
such that
$\Gamma_m\cap\partial D$
consists of at most one point.
In Section~\ref{s:KeyLemma},
we investigate
a $2$-angled disk whose interior
contains exactly three white vertices, and
we shall show a key lemma 
(Lemma~\ref{LemmaPartType322})
for Theorem~\ref{MainTheorem}.
In Section~\ref{s:IOC},
we review IO-Calculation
(a property of numbers of inward arcs of label $k$ and outward arcs of label $k$ in a closed domain $F$
with $\partial F\subset\Gamma_{k-1}\cup\Gamma_k\cup\Gamma_{k+1}$
for some label $k$).
In Section~\ref{s:Useful},
we introduce useful lemmata.
In Section~\ref{s:MainTheorem},
we prove Theorem~\ref{MainTheorem}.


\section{Preliminaries}
\label{s:Prel}

In this section, 
we introduce 
the definition of charts and its related words.

\begin{figure}[htb]
\begin{center}
\includegraphics{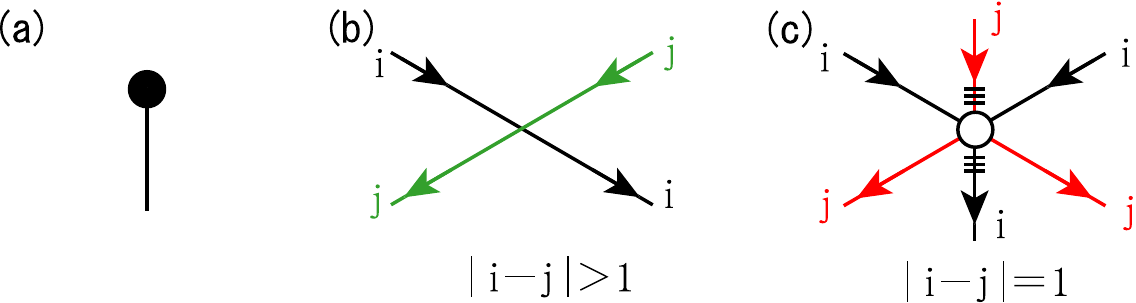}
\end{center}
\caption{ \label{fig01} (a) A black vertex. (b) A crossing. (c) A white vertex. 
Each arc with three transversal short arcs is a middle arc at the white vertex. }
\end{figure}

Let $n$ be a positive integer.
An $n$-{\it chart}  
(a braid chart of degree $n$ \cite{KnottedSurfaces}
or a surface braid chart of degree $n$ \cite{BraidBook}) 
is 
an oriented labeled graph in the interior of a disk,
which may be empty 
or
have closed edges without vertices
satisfying the following four conditions
(see Fig.~\ref{fig01}):
\begin{enumerate}
\item[(i)] 
Every vertex has degree $1$, $4$, or $6$.
\item[(ii)] 
The labels of edges are 
in $\{1,2,\dots,n-1\}$.
\item[(iii)]
In a small neighborhood of
each vertex of degree $6$,
there are six short arcs,
three consecutive arcs are
oriented inward 
and
the other three are outward,
and
these six are labeled $i$ and $i+1$
alternately for some $i$,
where the orientation and label of
each arc are inherited from
the edge containing the arc.
\item[(iv)]
For each vertex of degree $4$,
diagonal edges have the same label
and
are oriented coherently,
and the labels $i$ and $j$ of
the diagonals satisfy $|i-j|>1$.
\end{enumerate}
We call a vertex of degree $1$ a {\it black vertex},
a vertex of degree $4$ a {\it crossing}, and 
a vertex of degree $6$ a {\it white vertex}
respectively.
Among six short arcs
in a small neighborhood of
a white vertex,
a central arc of each three consecutive arcs
oriented inward (resp. outward) 
is called a   
{\it middle arc} at the white vertex
(see Fig.~\ref{fig01}(c)).
For each white vertex $v$, 
there are two middle arcs at $v$ 
in a small neighborhood of $v$.



\begin{figure}
\begin{center}
\includegraphics{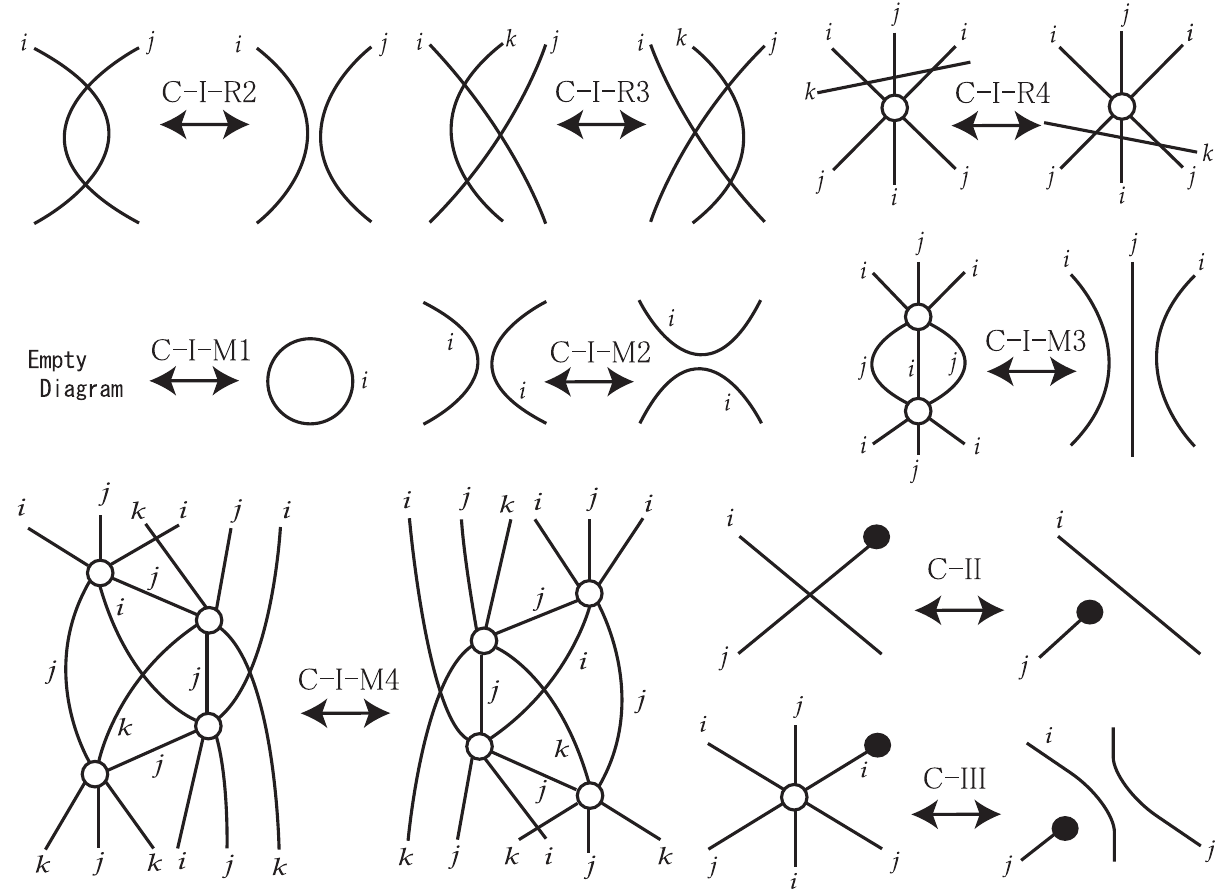}
\end{center}
\caption{ \label{fig02} For the C-III move, the edge containing the black vertex does not contain a middle arc at
a white vertex in the left figure. }
\end{figure}

Now {\it C-moves} are local modifications 
of charts as shown in Fig.~\ref{fig02}
(cf. \cite{KnottedSurfaces}, 
\cite{BraidBook} and \cite{Tanaka}).
Two charts are said to be {\it C-move equivalent}  if there exists
a finite sequence of C-moves 
which modifies one of the two charts 
to the other.

An edge in a chart is called 
a {\it free edge}
if it has
two black vertices.

For each chart $\Gamma$,
let $w(\Gamma)$ and $f(\Gamma)$ be the number of white vertices, and the number of free edges respectively.
The pair $(w(\Gamma), -f(\Gamma))$ is called a {\it complexity} of the chart (see \cite{BraidThree}).
A chart $\Gamma$ is called a {\it minimal chart} if its complexity is minimal among the charts C-move equivalent to the chart $\Gamma$ with respect to the lexicographic order of pairs of integers.

We showed the difference of a chart in a disk and in a 2-sphere (see \cite[Lemma 2.1]{ChartApp1}).
This lemma follows from that there exists a natural one-to-one correspondence between $\{$charts in $S^2\}/$C-moves and $\{$charts in $D^2\}/$C-moves, conjugations
(\cite[Chapter 23 and Chapter 25]{BraidBook}).
To make the argument simple, we assume that 
the charts lie on the 2-sphere instead of the disk.
\begin{assumption}
In this paper,
all charts are contained in the $2$-sphere $S^2$.
\end{assumption}
We have the special point in the 2-sphere $S^2$, called the point at infinity,
 denoted by $\infty$.
In this paper, all charts are contained in a disk such that the disk 
does not contain the point at infinity $\infty$.

An edge in a chart is called 
a {\it terminal edge}
if it has
a white vertex and a black vertex.

Let $\Gamma$ be a chart,
and $m$ a label of $\Gamma$. 
A {\it hoop} is a closed edge of $\Gamma$ without vertices 
(hence without crossings, neither).
A {\it ring} is a simple closed curve in $\Gamma_m$ containing a crossing but not containing any white vertices.
A hoop is said to be {\it simple} 
if one of the two complementary domains
of the hoop
does not contain any white vertices.

We can assume that
all minimal charts $\Gamma$
satisfy the following four conditions 
(see \cite{ChartApp1},\cite{ChartAppII},\cite{ChartAppIII},
\cite{StI}):

\begin{assumption}
\label{AssumeTerminal}
If an edge of $\Gamma$
contains a black vertex,
then the edge is a free edge 
or a terminal edge.
Moreover 
any terminal edge contains a middle arc.
\end{assumption}

\begin{assumption}
\label{NoSimpleHoop}
All free edges and simple hoops in $\Gamma$ 
are moved into a small neighborhood $U_\infty$ 
of the point at infinity $\infty$. 
Hence
we assume that 
$\Gamma$ does not contain free edges
nor simple hoops, 
otherwise mentioned. 
\end{assumption}

\begin{assumption}
\label{Ring}
Each complementary domain of
any ring and hoop must contain 
at least one white vertex. 
\end{assumption}

\begin{assumption}
The point at infinity $\infty$ is moved in any complementary domain of $\Gamma$.
\end{assumption}

In this paper
for a set $X$ in a space
we denote 
the interior of $X$,
the boundary of $X$ and
the closure of $X$
by Int$X$, $\partial X$
and $Cl(X)$
respectively.

\section{Connected components of $\Gamma_m$}
\label{s:ConnectedComponent}
In this section,
we investigate connected components of $\Gamma_m$
with at most three white vertices
for a minimal chart $\Gamma$.

\begin{figure}[hbt]
\centerline{\includegraphics{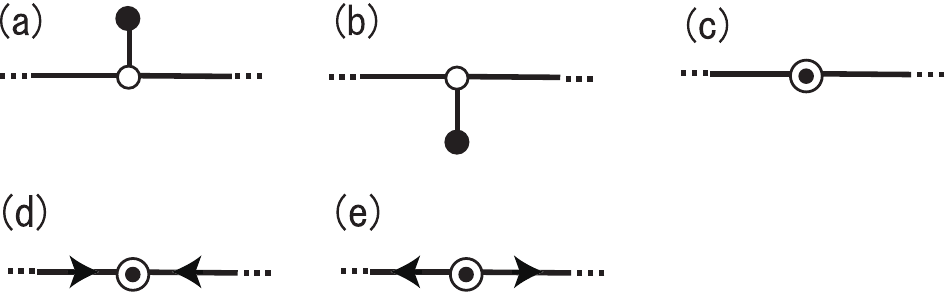}}
\caption{\label{fig03}
(a), (b) white vertices in terminal edges,
(c), (d), (e) BW-vertices.}
\end{figure}

 In our argument  we often construct a chart $\Gamma$. 
On the construction of a chart $\Gamma$, for a white vertex $w\in\Gamma_m$ for some label $m$,  
among the three edges of $\Gamma_m$ 
containing $w$, 
if one of the three edges is a terminal edge 
(see Fig.~\ref{fig03}(a) and (b)), 
then we remove the terminal edge and
put a black dot at the center of the white vertex  as shown in Fig.~\ref{fig03}(c).
Namely
Fig.~\ref{fig03}(c) means 
Fig.~\ref{fig03}(a) or 
Fig.~\ref{fig03}(b).
We call the vertex in Fig.~\ref{fig03}(c) 
a {\it BW-vertex}.

\begin{lemma}
\label{OriBWvertex}
In a minimal chart,
two edges containing the same BW-vertex
are oriented inward or outward at the BW-vertex
simultaneously
$($see Fig.~\ref{fig03}$($d$)$ and $($e$))$.
\end{lemma}

\begin{Proof}
By Assumption~\ref{AssumeTerminal},
each terminal edge of label $m$ contains a middle arc at a white vertex.
Thus the other two edges of label $m$ 
are oriented inward or outward at the BW-vertex
simultaneously.
Hence we have the result.
\end{Proof}

Let $\Gamma$ be a chart,
and $m$ a label of $\Gamma$. 
A {\it loop} is a simple closed curve in $\Gamma_m$ with exactly one white vertex
(possibly with crossings).

Let $X$ be a set in a chart $\Gamma$.
Let
 $$w(X)=\text{the number of white vertices in $X$.}$$

The following lemma is easily shown.
Thus we omit the proof.

\begin{figure}[htb]
\centerline{\includegraphics{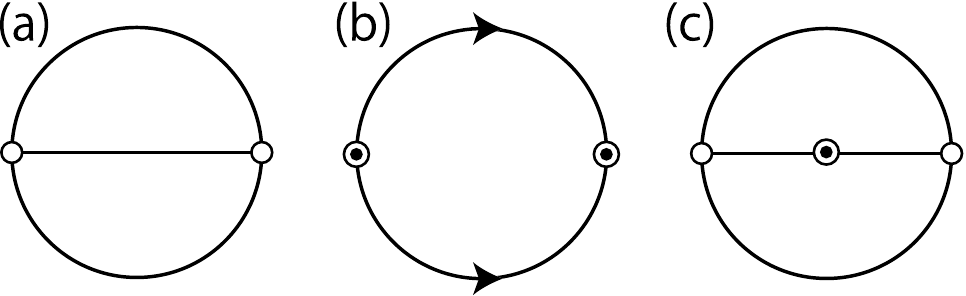}}
\caption{\label{fig04}
(a), (b) graphs with two white vertices, 
(c) a graph with three white vertices.}
\end{figure}

\begin{lemma}
\label{LemmaWithTerminal}
Let $\Gamma$ be a minimal chart,
and $m$ a label of $\Gamma$.
Let $G$ be a connected component of $\Gamma_m$.
Then we have the following.
\begin{enumerate}
\item[$(1)$] If $1\le w(G)$, then $2\le w(G)$.
\item[$(2)$] If $1\le w(G)\le 3$
and $G$ does not contain any loop, 
then $G$ is one of three graphs as shown 
in Fig.~\ref{fig04}.\hfill {$\square$}
\end{enumerate}
\end{lemma}

By Lemma~\ref{OriBWvertex},
we have an orientation of the graph 
as in Fig.~\ref{fig04}(b).

\section{$k$-angled disks}
\label{s:KAngledDisk}

In this section
we review properties of $k$-angled disks.

Let $\Gamma$ be a chart, $m$ a label of $\Gamma$, $D$ a disk with $\partial D\subset \Gamma_m$, 
and $k$ a positive integer.
If $\partial D$ contains exactly
$k$ white vertices, 
then $D$ is called 
{\it a $k$-angled disk of $\Gamma_m$}. 
Note that 
the boundary $\partial D$ may contain crossings.

Let $\Gamma$ be a chart, and
$m$ a label of $\Gamma$.
An edge of label $m$ is called a {\it feeler} of a $k$-angled disk $D$ of $\Gamma_m$
if the edge intersects $N-\partial D$
where $N$ is a regular neighborhood of $\partial D$ in $D$.

Let $\Gamma$ be a chart. 
Suppose that an object consists of 
some edges of $\Gamma$, arcs in edges of $\Gamma$ and arcs around white vertices.
Then the object is called {\it a pseudo chart}.

\begin{figure}[htb]
\centerline{\includegraphics{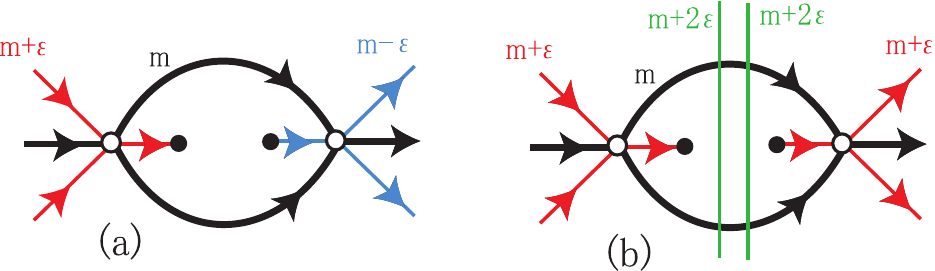}}
\caption{\label{fig05}
$m$ is a label,
and $\varepsilon\in\{+1,-1\}$.}
\end{figure}

\begin{lemma}
\label{Theorem2AngledDisk}
{\rm (\cite[Corollary 6.2]{ChartAppII})}
Let $\Gamma$ be a minimal chart.
Let $D$ be a $2$-angled disk of $\Gamma_m$ with at most one feeler.
If $w(\Gamma\cap{\rm Int}D)=0$,
then a regular neighborhood of $D$ contains one of two pseudo charts as shown in Fig.~\ref{fig05}.
\end{lemma}

\begin{figure}[hbt]
\centerline{\includegraphics{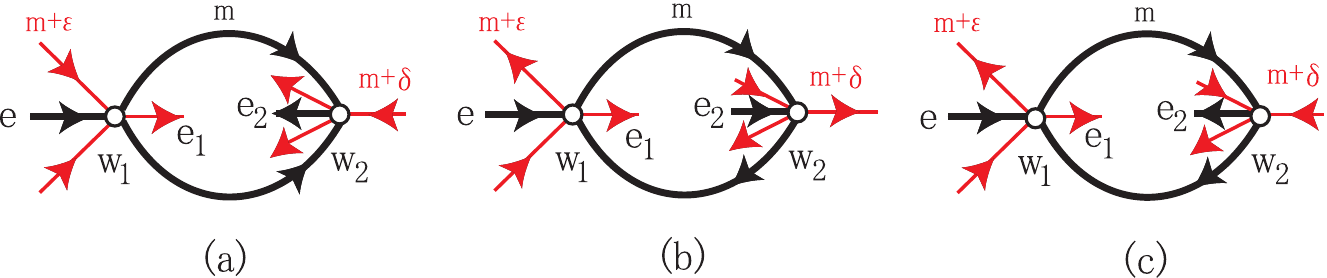}}
\caption{\label{fig06} 
The white vertex $w_1$ is in 
$\Gamma_m\cap\Gamma_{m+\varepsilon}$ and 
the white vertex $w_2$ is in 
$\Gamma_m\cap\Gamma_{m+\delta}$ 
where $\varepsilon,\delta\in\{+1,-1\}$.}
\end{figure}

Let $D$ be a 2-angled disk of $\Gamma_m$
 with exactly one feeler, 
 and $e$ an edge of label $m$ containing a white vertex $w_1$ in $\partial D$ but not contained in $D$.
If necessary we take the reflection of the chart $\Gamma$ or change the orientations of all of the edges, we have the following three 2-angled disks as shown in
Fig.~\ref{fig06}.

\begin{figure}[hbt]
\centerline{\includegraphics{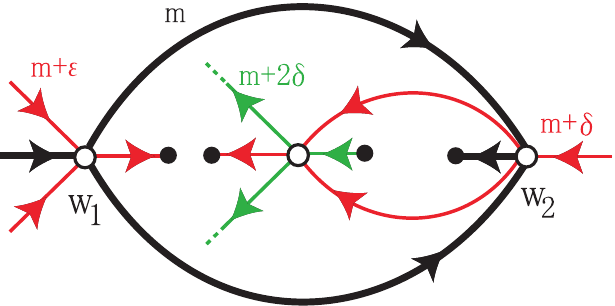}}
\caption{\label{fig07}
$m$ is a label, $\varepsilon,\delta\in\{+1,-1\}$.}
\end{figure}

\begin{lemma}
{\rm (\cite[Lemma 6.1]{ChartAppII})}
\label{Lemma2Cycle1-A}
Let $\Gamma$ be a minimal chart. 
Let $D$ be a $2$-angled disk of $\Gamma_m$ as shown in Fig.~\ref{fig06}$($a$)$.
Then $w(\Gamma\cap{\rm Int}D)\ge1$.
If  $w(\Gamma\cap{\rm Int}D)=1$, then 
a regular neighborhood of $D$ contains the pseudo chart as shown in Fig.~\ref{fig07}.
\end{lemma}

\begin{lemma}
{\rm (\cite[Lemma 6.2]{ChartAppII})}
\label{Lemma2Cycle1-B}
Let $\Gamma$ be a minimal chart. 
Let $D$ be a $2$-angled disk of $\Gamma_m$ as shown in Fig.~\ref{fig06}$($b$)$ or $($c$)$.
Then $w(\Gamma\cap{\rm Int}D)\ge3$. 
\end{lemma}

Let $\Gamma$ be a chart, 
$D$ a $k$-angled disk of $\Gamma_m$, and 
$G$ a pseudo chart with $\partial D\subset G$.
Let $r:D\to D$ be a reflection of $D$, and $G^*$ the pseudo chart obtained from $G$ by changing the orientations of all of the edges.
Then the set $\{G,G^*, r(G), r(G^*)\}$ 
is called the {\it RO-family of the pseudo chart $G$}.

By Lemma~\ref{Lemma2Cycle1-A} and Lemma~\ref{Lemma2Cycle1-B},
we have the following lemma.

\begin{lemma}
\label{OvalOneFeeler}
Let $\Gamma$ be a minimal chart,
and $m$ a label of $\Gamma$.
Let $D$ be a $2$-angled disk of $\Gamma_m$
with exactly one feeler.
Then we have the following:
\begin{enumerate}
\item[$(1)$] 
$w(\Gamma\cap {\rm Int}D)\ge1$. 
\item[$(2)$] 
If $w(\Gamma\cap {\rm Int}D)=1$,
then a regular neighborhood of $D$ contains 
one of the RO-family of
the pseudo chart as shown in 
Fig.~\ref{fig07}.
\end{enumerate}
\end{lemma}

\section{A disk with exactly two white vertices}
\label{s:DiskTwo}
In this section
for a minimal chart $\Gamma$
we investigate
a disk $D$ with exactly two white vertices 
of $\Gamma_m$
such that
$\Gamma_m\cap\partial D$
consists of at most one point.

Let $\Gamma$ be a chart, 
and $m$ a label of $\Gamma$.
Let $L$ be the closure of a connected component 
of the set obtained by taking out 
all the white vertices from $\Gamma_m$.
If $L$ contains at least one white vertex
but does not contain any black vertex,
then $L$ is called an {\it internal edge of label $m$}.
Note that an internal edge may contain a crossing of $\Gamma$.

\begin{figure}[htb]
\centerline{\includegraphics{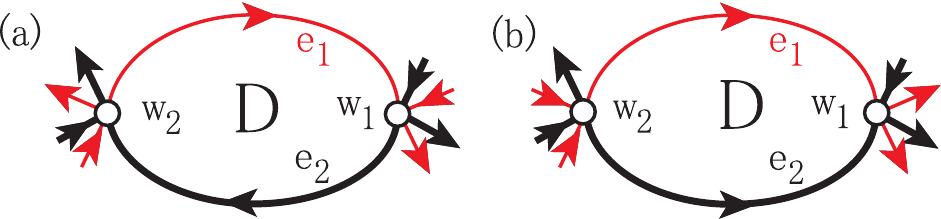}}
\caption{\label{fig08}
Lenses.}
\end{figure}

Let $\Gamma$ be a chart. 
Let $D$ be a disk 
such that 
\begin{enumerate}
\item[(1)] $\partial D$ consists of an internal edge $e_1$ of label $m$ and an internal edge $e_2$ of label ${m+1}$, and 
\item[(2)] any edge containing a white vertex in $e_1$ does not intersect the open disk Int$D$.
\end{enumerate}
Note that $\partial D$ may contain crossings.
Let $w_1$ and $w_2$ be the white vertices in $e_1$. 
If the disk $D$ satisfies one of the following conditions, then $D$ is called  {\it a lens of type $(m,m+1)$}
(see Fig.~\ref{fig08}):
\begin{enumerate}
	\item[(i)] Neither $e_1$ nor $e_2$ contains a middle arc. 
	\item[(ii)] One of the two edges $e_1$ and $e_2$ contains middle arcs at both white vertices $w_1$ and $w_2$ simultaneously.
\end{enumerate}

\begin{lemma}{\rm (\cite[Theorem 1.1]{ChartApp1}
and \cite[Corollary 1.1]{ChartAppII})}
\label{LemmaNoLens}
\begin{enumerate}
\item[$(1)$] 
There exist at least three white
vertices in the interior of a lens
for any minimal chart.
\item[$(2)$] There is no lens in any minimal chart with 
at most seven white vertices.
\end{enumerate}
\end{lemma}

\begin{figure}[htb]
\centerline{\includegraphics{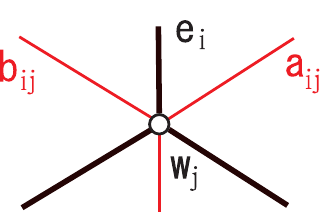}}
\caption{\label{fig09}
The three edges $a_{ij},e_i,b_{ij}$ are consecutive edges around the white vertex $w_j$.}
\end{figure}

In our argument,
we often need a name for an unnamed edge by using a given edge and a given white vertex.
For the convenience,
we use the following naming:
Let $e',e_i,e''$ be three consecutive edges containing  a white vertex $w_j$. Here, 
the two edges $e'$ and $e''$ are unnamed edges. 
There are six arcs in a neighborhood $U$ of the white vertex $w_j$. 
If the three arcs $e'\cap U$, $e_i \cap U$, $e'' \cap U$ lie anticlockwise around the white vertex $w_j$ in this order, 
then $e'$ and $e''$ are denoted by $a_{ij}$ and $b_{ij}$ 
respectively (see Fig.~\ref{fig09}).
There is a possibility $a_{ij}=b_{ij}$ if they are contained in a loop.

\begin{lemma}
\label{OvalType4}
Let $\Gamma$ be a minimal chart,
and $m$ a label of $\Gamma$.
Let $D$ be a $2$-angled disk of $\Gamma_m$.
Then we have the following:
\begin{enumerate}
\item[$(1)$] 
If there exist two feelers of $D$
each of which is a terminal edge,
then $w(\Gamma\cap {\rm Int}D)\ge1$. 
\item[$(2)$] 
If $w(\Gamma\cap {\rm Int}D)=0$,
then $D$ has at most one feeler.
\end{enumerate}
\end{lemma}

\begin{figure}[hbt]
\centerline{\includegraphics{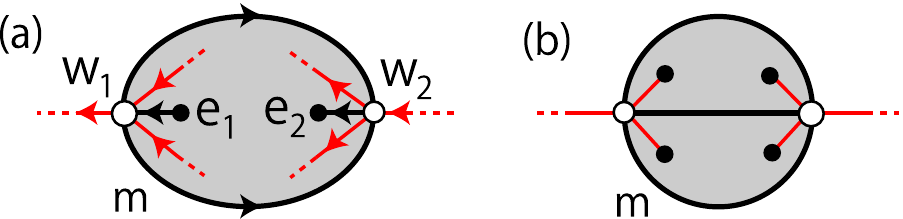}}
\caption{\label{fig10}
The gray region is the disk $D$.}
\end{figure}

\begin{Proof}
We show Statement (1).
Let $e_1,e_2$ be feelers of $D$.
By the condition of this lemma,
the edges $e_1,e_2$ are terminal edges
(see Fig.~\ref{fig10}(a)).
Let $w_1,w_2$ be the white vertices in $e_1,e_2$
respectively.
Without loss of generality
we can assume that 
$e_1$ is oriented inward at $w_1$.
Since the terminal edge $e_1$ contains a middle arc at $w_1$
by Assumption~\ref{AssumeTerminal},
we have orientation of the other edges
as shown in Fig.~\ref{fig10}(a).

Suppose $w(\Gamma\cap{\rm Int}D)=0$.
Let $a_{11},b_{11}$ be internal edges 
(possibly terminal edges)
of label $m+\varepsilon$
in $D$ with $w_1\in a_{11}\cap b_{11}$,
here $\varepsilon\in\{+1,-1\}$.
Since neither $a_{11}$ nor $b_{11}$
contains a middle arc at $w_1$,
by Assumption~\ref{AssumeTerminal}
neither $a_{11}$ nor $b_{11}$
is a terminal edge.
Hence both of $a_{11}$ and $b_{11}$
contain the white vertex $w_2$.
Thus $w(\Gamma\cap{\rm Int}D)=0$ implies that
there are two lenses of 
type $(m,m+\varepsilon)$ in $D$
whose interiors do not contain 
any white vertices.
This contradicts Lemma~\ref{LemmaNoLens}(1).
Hence $w(\Gamma\cap {\rm Int}D)\ge1$.

We show Statement (2).
Suppose that $D$ has two feelers $e_1,e_2$.
Since $w(\Gamma\cap{\rm Int}D)=0$,
we have that 
$e_1=e_2$ or both of $e_1,e_2$ are terminal edges.

If $e_1=e_2$,
i.e. the set $e_1\cup \partial D$ is 
a connected component of 
$\Gamma_m$ as in Fig.~\ref{fig04}(a),
then the disk $D$ separates into two $2$-angled disks 
of $\Gamma_m$ without feelers.
By Lemma~\ref{Theorem2AngledDisk},
each $2$-angled disk contains one of 
two pseudo charts as in Fig.~\ref{fig05}
(see Fig.~\ref{fig10}(b)).
Thus at each white vertex in $D$
there exist at least two terminal edges 
each of whose label is different from $m$.
Namely, at each white vertex in $D$
there exist at least two terminal edges
of the same label.
One of the two terminal edges does not contain
a middle arc.
This contradicts Assumption~\ref{AssumeTerminal}.

If
both of $e_1,e_2$ are terminal edges,
then by Lemma~\ref{OvalType4}(1)
we have $w(\Gamma\cap{\rm Int}D)\ge1$.
This is a contradiction.
Hence 
$D$ has at most one feeler.
\end{Proof}

\begin{figure}[htb]
\centerline{\includegraphics{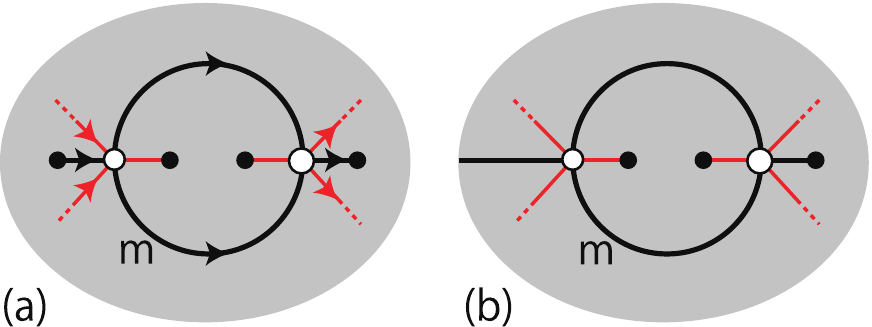}}
\caption{\label{fig11}
The gray region is the disk $D$,
and $m$ is a label.}
\end{figure}

\begin{lemma}
\label{LemmaTwoWhiteVertices}
Let $\Gamma$ be a minimal chart without loops, 
and
$m$ a label of $\Gamma$.
Let $D$ be a disk such that
if an edge intersects 
$\partial D$,
then the edge intersects $\partial D$
transversely.
Suppose that
$D$ contains exactly two white vertices $w_1,w_2$ and 
$\Gamma_m\cap \partial D$ is at most one point.
Then we have the following:
\begin{enumerate}
\item[$(1)$] 
If $w_1,w_2\in \Gamma_{m}$,
then the disk $D$ contains one of the 
two pseudo charts
as shown in 
Fig.~\ref{fig11}.
\item[$(2)$]
If there exists a number $\varepsilon\in\{+1,-1\}$
such that 
$w_1,w_2\in \Gamma_m\cap\Gamma_{m+\varepsilon}$
and 
if $S^2-D$ does not contain any white vertices in 
$\Gamma_{m+\varepsilon}$,
then there exist two lenses of 
type $(m,m+\varepsilon)$.
\end{enumerate}
\end{lemma}


\begin{Proof}
We shall show Statement (1).
Let $G$ be a connected component of 
$\Gamma_m\cap D$
with $w(G)\ge1$.

{\bf Case (i).} $G\cap\partial D=\emptyset$.

Since $D$ contains exactly two white vertices,
we have $w(G)\le2$.
By Lemma~\ref{LemmaWithTerminal}(2),
the graph $G$ is one of the two graphs as in
Fig.~\ref{fig04}(a),(b).

If $G$ is the graph as in 
Fig.~\ref{fig04}(a),
then the disk $D$ contains two $2$-angled disk 
of $\Gamma_m$ without feelers.
By Lemma~\ref{Theorem2AngledDisk},
each $2$-angled disk contains one of 
two pseudo charts as in Fig.~\ref{fig05}
(see Fig.~\ref{fig10}(b)).
By a similar way to the proof of 
Lemma~\ref{OvalType4}(2),
we have a contradiction.

If $G$ is the graph as in 
Fig.~\ref{fig04}(b),
then there is a $2$-angled disk $D'$ of $\Gamma_m$
in $D$. 
Now $w(\Gamma\cap D)=2$
implies that ${\rm Int}D'$ does not contain 
any white vertices.
Thus by Lemma~\ref{OvalType4}(2),
the 2-angled disk $D'$ has at most one feeler.
Hence by Lemma~\ref{Theorem2AngledDisk}
the 2-angled disk $D'$ contains one of 
 two pseudo charts as in Fig.~\ref{fig05}.
Thus $D$ contains the pseudo chart as in
Fig.~\ref{fig11}(a).

{\bf Case (ii).} $G\cap\partial D\not=\emptyset$.

Let $e$ be the internal edge 
(possibly terminal edge) containing 
the point $G\cap\partial D$,
and $v$ the endpoint of $e$ with $v\not\in D$.

If $v$ is a black vertex,
then we can show that 
the disk $D$ contains the pseudo chart as in
Fig.~\ref{fig11}(a)
 by a similar way to Case (i).
If $v$ is a white vertex,
then the disk $D$ contains a $2$-angled disk $D'$
of $\Gamma_m$,
because
$\Gamma$ has no loop.
Since $e$ is not a feeler of $D'$,
the disk $D'$ has at most one feeler.
Hence by the similar way as the one of Case (i),
we can show that the 2-angled disk $D'$ contains 
one of 
 two pseudo charts as in Fig.~\ref{fig05}.
Thus $D$ contains the pseudo chart as in
Fig.~\ref{fig11}(b).

We shall show Statement (2).
By Statement (1),
the disk $D$ contains one of 
the two pseudo charts as in
Fig.~\ref{fig11}.

If $D$ contains the pseudo chart as in
Fig.~\ref{fig11}(a),
then the disk $D$ contains a $2$-angled disk $D'$ 
of $\Gamma_m$ without feelers.
Thus $Cl(S^2-D')$ is a $2$-angled disk
with two feelers
each of which is a terminal edge.
Since 
$S^2-D$ does not contain any white vertices 
in $\Gamma_{m+\varepsilon}$, 
the set $S^2-D'$ does not contain 
any white vertex of $\Gamma_{m+\varepsilon}$.
By a similar way to 
the proof of Lemma~\ref{OvalType4}(1),
we can show that 
there exist two lenses of type $(m,m+\varepsilon)$.

If $D$ contains the pseudo chart as in
Fig.~\ref{fig11}(b),
then similarly
we can show that 
there exist two lenses of type $(m,m+\varepsilon)$.
\end{Proof}

\section{2-angled disks of $\Gamma_k$
whose interiors contain exactly three white vertices}
\label{s:KeyLemma}

In this section,
we investigate a 2-angled disk of $\Gamma_k$
whose interior contains exactly three white vertices,
and 
we shall show a key lemma
(Lemma~\ref{LemmaPartType322})
for Theorem~\ref{MainTheorem}.

\begin{figure}[htb]
\centerline{\includegraphics{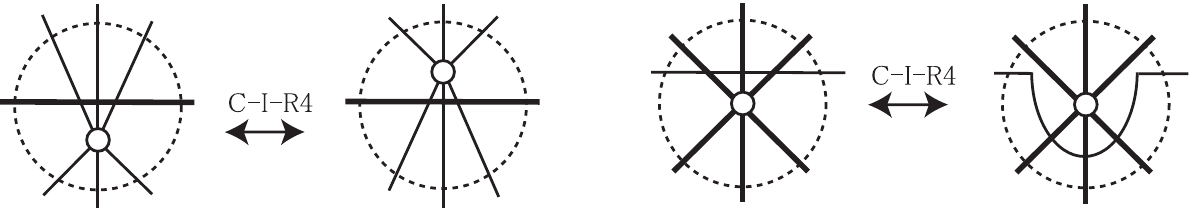}}
\caption{\label{fig12} 
C-moves keeping thicken figures fixed.}
\end{figure}

Let $\Gamma$ and $\Gamma^\prime $ be C-move equivalent charts. 
Suppose that a pseudo chart $X$ of $\Gamma$ is also a pseudo chart of $\Gamma^\prime$. 
Then we say that 
$\Gamma$ is modified to $\Gamma^\prime$ by {\it C-moves keeping $X$ fixed}.
In Fig.~\ref{fig12},
we give examples of C-moves keeping pseudo charts  fixed.

Let $\alpha$ be a simple arc,
and $p,q$ the endpoints of $\alpha$.
We denote 
$\partial \alpha=\{p,q\}$ and ${\rm Int}\alpha=\alpha-\partial \alpha$.

Let $\Gamma$ be a chart, and $D$ a disk.
Let $\alpha$ be a simple arc in $\partial D$,
and $\gamma$ a simple arc in an internal edge of label $k$.
The simple arc $\gamma$ is called
a {\it {$(D,\alpha)$-arc}} of label $k$
provided that 
$\partial \gamma \subset $Int$\alpha$
and
Int$\gamma\subset $Int$D$. 
If there is no $(D,\alpha)$-arc in $\Gamma$,
then the chart $\Gamma$ is said to be
$(D,\alpha)$-{\it arc free}.

The following lemma will be used in the
proof of Lemma~\ref{LemmaPartType322}.

\begin{figure}
\centerline{\includegraphics{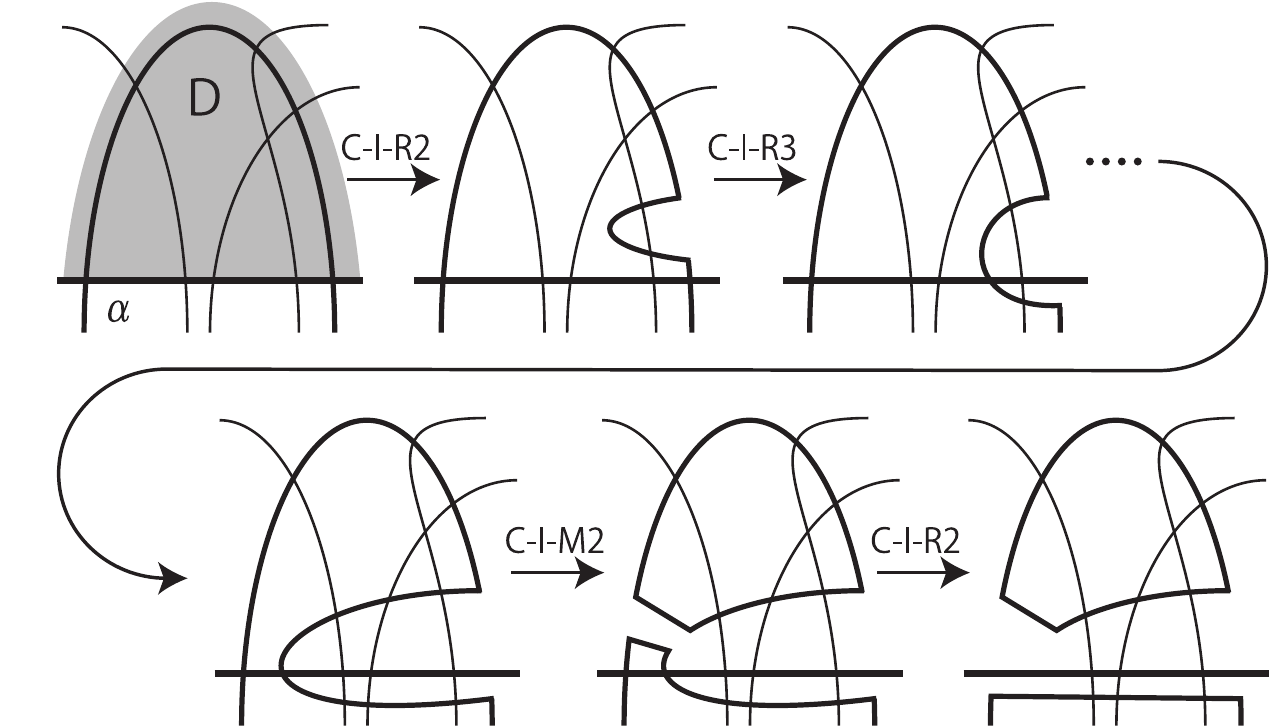}}
\caption{\label{fig13}
The gray region is the disk $D$.}
\end{figure}

\begin{lemma}
$($New Disk Lemma$)$
{\em (\cite[Lemma 7.1]{Chart33},
cf. \cite[Lemma 3.2]{ChartApp1})} 
\label{NewDiskLemma}
Let $\Gamma$ be a chart and
$D$ a disk 
whose interior does not contain 
a white vertex nor a black vertex of $\Gamma$.
Let $\alpha$ be a simple arc in $\partial D$ 
such that ${\rm Int}\alpha$ does not contain 
a white vertex nor a black vertex of $\Gamma$.
Let $V$ be a regular neighborhood of $\alpha$. 
Suppose that 
the arc $\alpha$ is contained in 
an internal edge of some label $k$ of $\Gamma$.
Then by applying C-I-M2 moves, C-I-R2 moves, 
and C-I-R3 moves in $V$, 
there exists 
a $(D,\alpha)$-arc free chart $\Gamma'$ 
obtained from the chart $\Gamma$ 
keeping $\alpha$ fixed 
$($see Fig.~\ref{fig13}$)$.
\end{lemma}

The following lemma will be used in the
proof of Lemma~\ref{Lemma2PartType322}.

\begin{lemma}
$(${\rm \cite[Lemma 6.1]{ChartAppIII}}$)$
\label{ConditionRing} 
Let $\Gamma$ be a minimal chart.
Let $C$ be a ring or a non simple hoop, 
and $D$ a disk with $\partial D=C$.
If $w(\Gamma\cap D)=1$, then $\Gamma$ is C-move equivalent to the minimal chart $Cl(\Gamma-C)$. 
\end{lemma}

\begin{figure}[hbt]
\centerline{\includegraphics{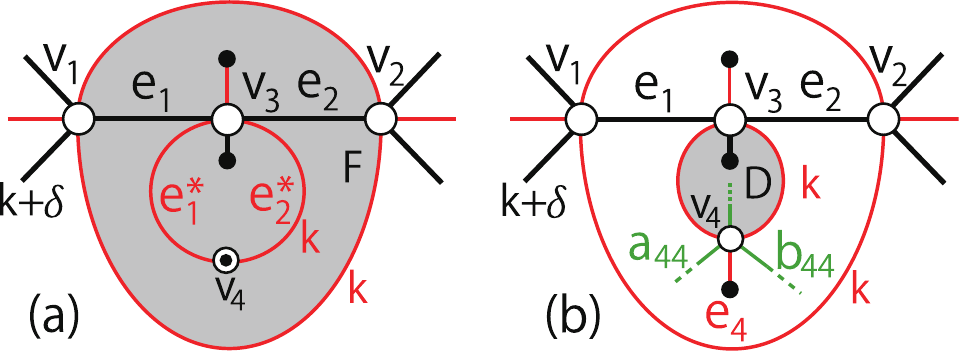}}
\caption{\label{fig14}
$k$ is the label, and $\delta\in\{+1,-1\}$.
(a) The gray region is the disk $F$.
(b) The gray region is the disk $D$.
}
\end{figure}

From now on throughout this section, we may assume
that 
\begin{enumerate}
\item[(i)] 
$\Gamma$ is a minimal chart, 
\item[(ii)] $F$ is 
 a $2$-angled disk of $\Gamma_k$
without feelers
with $w(\Gamma\cap{\rm Int}F)=3$
such that
a regular neighborhood of $F$ contains 
the pseudo chart 
as shown in Fig.~\ref{fig14}$($a$)$
where 
\begin{enumerate}
\item[{\rm (a)}] 
$v_1,v_2,v_3,v_4$ are white vertices in $F$ with 
$v_1,v_2\in\partial F$, 
$v_1,v_2,v_3\in\Gamma_k\cap\Gamma_{k+\delta}$,
and $v_4\in\Gamma_k\cap\Gamma_{k-\delta}$
here $\delta\in\{+1,-1\}$,
\end{enumerate}
\item[(iii)]
$v_5$ is the white vertex in
${\rm Int}F$ different from $v_3,v_4$ with
$v_5\in\Gamma_{k-\delta}\cap\Gamma_{k-2\delta}$.
\end{enumerate}

Let $\alpha$ be a simple arc,  
and $p,q$ points in $\alpha$. 
We denote by $\alpha[p,q]$ the subarc of $\alpha$ whose endpoints are $p$ and $q$. 

Let $D$ be a compact surface.
A simple arc $\alpha$ in $D$
is a {\it proper arc} of $D$
if $\alpha\cap\partial D=\partial \alpha$.

\begin{figure}[hbt]
\centerline{\includegraphics{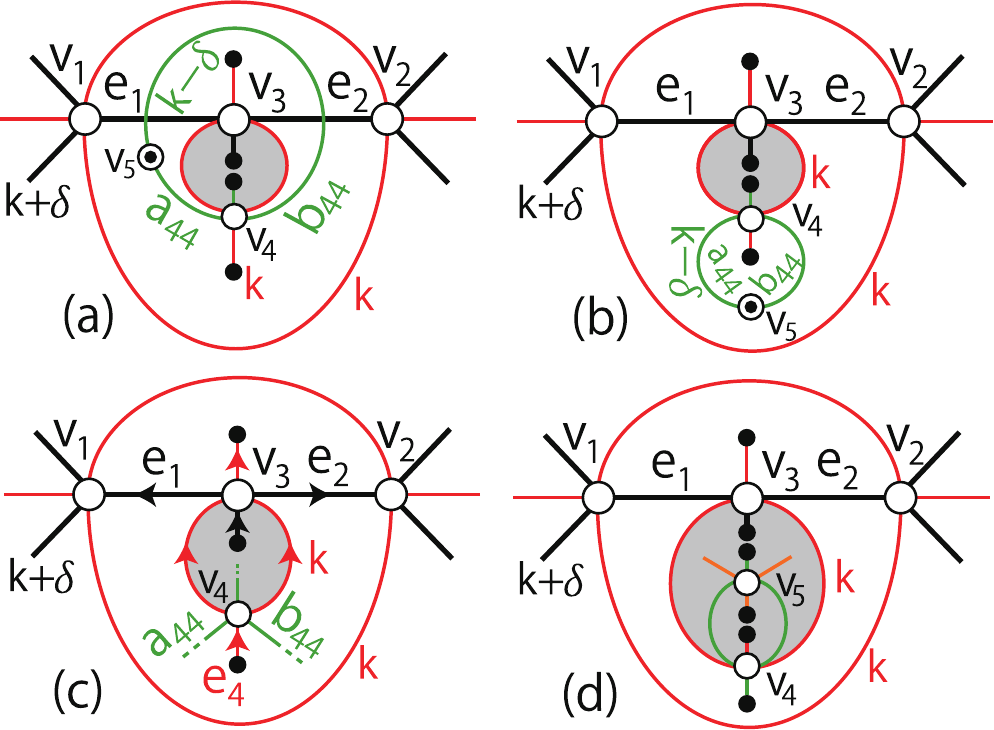}}
\caption{\label{fig15}
The gray regions are the disk $D$,
$k$ is the label, and $\delta\in\{+1,-1\}$.
(a)  $(e_1\cup e_2)\cap(a_{44}\cup b_{44})=$ two points. 
(b) $(e_1\cup e_2)\cap(a_{44}\cup b_{44})=\emptyset$.
(c) The terminal edge $e_4$ of label $k$ is oriented inward at $v_4$.
(d) The terminal edge $e_4$ of label $k$ is contained in $D$.
}
\end{figure}

\begin{lemma}
\label{LemmaPartType322}
Let $\Gamma,F,v_1,\cdots,v_5$ be as above. 
Let $D$ be the $2$-angled disk of $\Gamma_{k}$
in $F$ with $v_3,v_4\in\partial D$,
and $e_4$ the terminal edge 
at $v_4$ of label $k$.
If $e_4\not\subset D$ 
$($see Fig.~\ref{fig14}$($b$))$,
then $\Gamma$ can be modified to a minimal chart 
containing one of two pseudo charts as shown in
Fig.~\ref{fig15}$($a$)$ and $($b$)$
by C-moves in $F$ 
keeping $\Gamma_k\cup\Gamma_{k+\delta}$ fixed.
\end{lemma}


\begin{Proof}
We can assume that
$e_4$ is oriented inward at $v_4$.
Since $e_4$ contains a middle arc at $v_4$
by Assumption~\ref{AssumeTerminal},
we have orientation of edges as shown in
Fig.~\ref{fig15}(c).

Let $a_{44},b_{44}$ be the internal edges 
(possibly terminal edges) of 
label $k-\delta$
oriented inward at $v_4$
such that $a_{44},e_4,b_{44}$ lie anticlockwise
around the vertex $v_4$ in this order.
Since neither $a_{44}$ nor $b_{44}$
contains a middle arc at $v_4$,
neither $a_{44}$ nor $b_{44}$ is
a terminal edge
by Assumption~\ref{AssumeTerminal}.
Hence
${\rm Int}F\ni v_3,v_4,v_5$
and $w(\Gamma\cap {\rm Int}F)=3$ imply
 $a_{44}\cap b_{44}\ni v_5$.

{\bf Claim~1.}
We can assume that 
the edge $a_{44}$ does not contain any crossings 
by applying C-moves in $F$
keeping 
$\Gamma_k\cup\Gamma_{k+\delta}$ fixed.

{\it Proof of Claim~$1$.}
We can assume $\delta=+1$
(for the case $\delta=-1$,
 we can show the claim similarly).

Since the edge $a_{44}$ is of 
label $k-1$,
the set ${\rm Int}(a_{44})$ does not intersect
edges of label $k,k-1,k-2$.
Let $x$ be a point in $a_{44}$
such that $a_{44}[x,v_4]$ 
does not contain any crossings.
By C-I-R2 moves and C-I-R3 moves
keeping $a_{44}\cup(\cup_{j\ge k+1}\Gamma_j)$
fixed,
we can move each crossing in 
$a_{44}\cap(\cup_{i\le k-3}\Gamma_i)$
into $a_{44}[x,v_4]$
one by one
(see Fig.~\ref{fig16}(a),(b),(c)).
Here we use the notation $\Gamma$
for the modified chart.
By C-I-R2 moves and C-I-R4 moves,
we can move out all the crossings on 
$a_{44}[x,v_4]$
by passing through the vertex $v_4$
(see Fig.~\ref{fig16}(d),(e)).
Here we use the notation $\Gamma$
for the modified chart.
Thus each crossing in $a_{44}$
is contained in 
$a_{44}[v_5,x]\cap(\cup_{j\ge k+1}\Gamma_j)$.

Finally we shorten the edge $a_{44}$
to $a_{44}[x,v_4]$
by C-I-R2 moves and C-I-R4 moves.
We abuse the notation $a_{44}$ 
for the shortened edge.
Then we obtain the edge $a_{44}$
without crossings
(see Fig.~\ref{fig16}(f)).
Therefore Claim~$1$ holds.

\begin{figure}
\centerline{\includegraphics{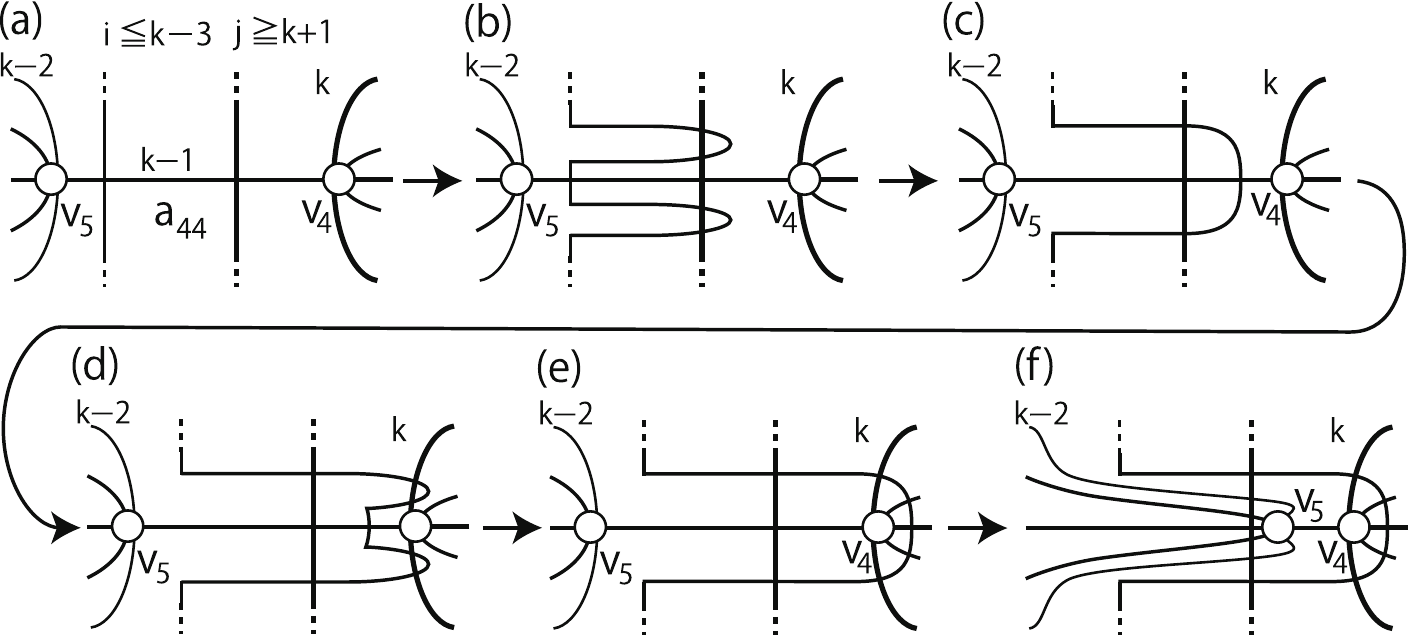}}
\caption{\label{fig16}
C-moves keeping thicken figures fixed.}
\end{figure}

Now the edge $b_{44}$ elongates as 
the edge $a_{44}$ shortens.
We also abuse the notation $b_{44}$ 
for the elongated edge.

The disk $D$ is the 2-angled disk of $\Gamma_k$
with $D\subset F$ and
$v_3,v_4\in\partial D$
(see Fig.~\ref{fig15}(c)).
Let $N(a_{44})$ be a regular neighborhood 
of $a_{44}$ in $F$.
Let ${\Bbb S}$ be the set of all minimal charts
each of which is modified from $\Gamma$
by C-moves in $F-D\cup N(a_{44})$
keeping $\Gamma_k\cup\Gamma_{k+\delta}$
fixed.
We can assume that $\Gamma$ is a minimal chart in ${\Bbb S}$
with
\begin{enumerate}
\item[(1)] $|(e_1\cup e_2)\cap \Gamma_{k-\delta}|={\rm min}\{\ 
|(e_1\cup e_2)\cap \Gamma_{k-\delta}' | 
\ : \ \Gamma'\in{\Bbb S}\}$
\end{enumerate}
where $|X|$ is the number of points in a set $X$.
Let $E$ be a disk in $F$ bounded by $a_{44}\cup b_{44}$.
Since $D\not\supset a_{44}$ and
$\partial D\cap \partial E=v_4$,
there are two cases:
(i) $E\supset D$,
(ii) $E\cap D=v_4$.

{\bf Case (i).}
Since $v_3\in {\rm Int }E$
and $v_1,v_2\not\in E$,
we have $|e_1\cap \partial E|\ge1$ and 
$|e_2\cap \partial E|\ge1$.
Thus 
$$|(e_1\cup e_2)\cap \Gamma_{k-\delta}|\ge
|(e_1\cup e_2)\cap (a_{44}\cup b_{44})|=
|(e_1\cup e_2)\cap \partial E|\ge2.$$

Now we show 
$|(e_1\cup e_2)\cap (a_{44}\cup b_{44})|=2$.
Suppose that 
$|(e_1\cup e_2)\cap  (a_{44}\cup b_{44}) |>2$.
Since $\partial E=a_{44}\cup b_{44}$,
we have $|(e_1\cup e_2)\cap \partial E |>2$.
The vertices $v_1,v_2$ are 
the endpoints of the arc
$e_1\cup e_2$.
Since $v_1$ and $v_2$ are outside the disk $E$,
the set $(e_1\cup e_2)\cap E$
consists of proper arcs of $E$.
Let $G$ be the proper arc containing $v_3$,
and $\alpha$ a proper arc different from $G$
(see Fig.~\ref{fig17}(a)).
The arc $\alpha$ divides the disk $E$
into two disks.
One of the two disks
does not intersect $G$,
say $E'$.

{\bf Claim 2.} $w(\Gamma\cap E')=0$.

{\it Proof of Claim~$2$.}
By Claim $1$, 
we have $(e_1\cup e_2)\cap a_{44}=\emptyset$.
Hence $(e_1\cup e_2)\cap (D\cup a_{44})=v_3$.
Thus $\alpha\cap (D\cup a_{44}\cup G)=\emptyset$.
Since $D\cup a_{44}\cup G$ is connected
and since $E'$ does not intersect $G$,
we have $E'\cap (D\cup a_{44}\cup G)=\emptyset$.
Thus $v_3,v_4,v_5\in D\cup a_{44}\cup G$
implies $w(\Gamma\cap E')=0$.
Hence Claim~$2$ holds.

Let $N(E')$ be a regular neighborhood of $E'$ in $F$.
Then the disk $N(E')$ contains a proper arc $\widetilde\alpha$ of label $k+\delta$ containing $\alpha$.
Let $\widetilde E$ be the disk divided by 
$\widetilde\alpha$ from $N(E')$
with $\widetilde E\supset E'$.
Applying New Disk Lemma (Lemma~\ref{NewDiskLemma})
for the disk $\widetilde E$,
we obtain a 
$(\widetilde E,\widetilde\alpha)$-arc free 
minimal chart $\widetilde \Gamma$
(see Fig.~\ref{fig17}(b)).
Thus $$|(e_1\cup e_2)\cap \widetilde\Gamma_{k-\delta} |<|(e_1\cup e_2)\cap \Gamma_{k-\delta}|.$$
This contradicts Condition (1).
Hence
$|(e_1\cup e_2)\cap (a_{44}\cup b_{44}) |=2$.
Hence
a regular neighborhood of 
$F$ contains the pseudo chart
as in Fig.~\ref{fig15}(a).

{\bf Case (ii).}
Similarly
we can show 
$(e_1\cup e_2)\cap (a_{44}\cup b_{44}) =\emptyset$
by modifying the chart $\Gamma$ 
by C-moves.
Thus
a regular neighborhood of $F$ 
contains the pseudo chart
as in Fig.~\ref{fig15}(b).
\end{Proof}

\begin{figure}
\centerline{\includegraphics{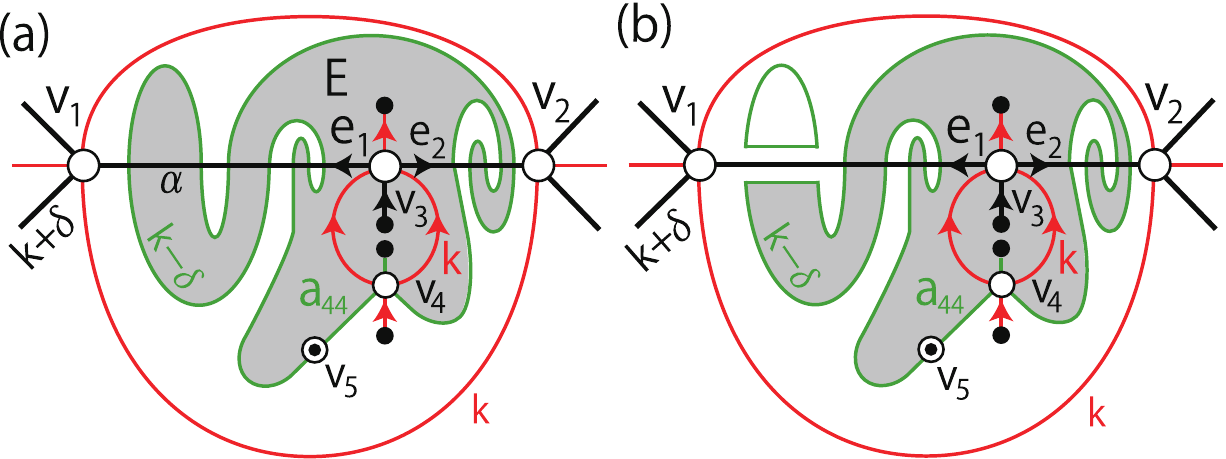}}
\caption{\label{fig17}
(a) The gray region is the disk $E$.
(b) $(\widetilde E,\widetilde\alpha)$-arc free 
minimal chart $\widetilde \Gamma$.}
\end{figure}

\begin{lemma}
\label{Lemma3PartType322}
Let $\Gamma,F,v_1,\cdots,v_5$ be as above.
If $\Gamma$ contains the pseudo chart as shown in
Fig.~\ref{fig14}$($a$)$,
then
the chart $\Gamma$ can be modified to a 
minimal chart containing 
one of three pseudo charts as shown in
Fig.~\ref{fig15}$($a$),($b$),($d$)$
by C-moves in $F$ keeping 
$\Gamma_k\cup\Gamma_{k+\delta}$ fixed.
\end{lemma}


\begin{Proof}
Let $D$ be the 2-angled disk of $\Gamma_k$
with $D\subset {\rm Int}F$ and
$v_3,v_4\in \partial D$.
Since $w(\Gamma\cap {\rm Int }F)=3$,
we have
$$3=w(\Gamma\cap {\rm Int }F)\ge
w(\Gamma\cap \partial D)+w(\Gamma\cap {\rm Int }D)= 2+w(\Gamma\cap {\rm Int }D).$$ 
Hence 
\begin{enumerate}
\item[(1)] $w(\Gamma\cap {\rm Int }D)\le 1$.
\end{enumerate}
Let $e_4$ be the terminal edge of label $k$
at $v_4$.
There are two cases: 
$e_4\subset D$ or $e_4\not\subset D$.

If $e_4\not\subset D$,
then by Lemma~\ref{LemmaPartType322}
the chart $\Gamma$ can be modified to a minimal chart
containing one of the two pseudo charts
as in Fig.~\ref{fig15}(a),(b)
by C-moves in $F$ keeping $\Gamma_k\cup\Gamma_{k+\delta}$ fixed.

If $e_4\subset D$,
then by (1) and  Lemma~\ref{OvalOneFeeler}(1)
we have $w(\Gamma\cap {\rm Int }D)= 1$.
Thus by Lemma~\ref{OvalOneFeeler}(2)
a regular neighborhood of $D$
contains the pseudo chart
as in Fig.~\ref{fig07}.
Hence $\Gamma$ contains 
the pseudo chart as shown in
Fig.~\ref{fig15}(d).
\end{Proof}

\begin{lemma}
\label{Lemma2PartType322}
Let $\Gamma,F,v_1,\cdots,v_5$ be as above.
Let $G$ be the union
of all the internal edges of label 
$k-\delta, k, k+\delta$ in $F$.
Suppose that a regular neighborhood of $F$ 
contains one of the three pseudo charts 
as shown in Fig.~\ref{fig15}$($a$),($b$),($d$)$.
Then
the chart $\Gamma$ can be modified to 
a minimal chart by C-moves keeping $G$
fixed so that there is no ring of label 
$k-\delta, k, k+\delta$ in $F$.
\end{lemma}

\begin{Proof}
By the condition of this lemma, 
\begin{enumerate}
\item[(1)] 
a regular neighborhood of $F$ contains one of 
three pseudo charts as in
Fig.~\ref{fig15}(a),(b),(d).
\end{enumerate}
Hence a regular neighborhood of $F$
contains the pseudo chart as in 
Fig.~\ref{fig14}(a).
We use the notations as in 
Fig.~\ref{fig14}(a)
where
$e_1,e_2$ are internal edges of 
label $k+\delta$, and
$e_1^*,e_2^*$ are internal edges of label $k$.
Let $e_4$ be the terminal edge of label $k$ at $v_4$.
Let $a_{44},b_{44}$ be internal edges of 
label $k-\delta$
containing $v_4$
such that $a_{44},e_4,b_{44}$
lie anticlockwise around $v_4$ 
in this order.
By (1),
the six edges 
$e_1,e_2,e_1^*,e_2^*,a_{44},b_{44}$
and two internal edges in $\partial F$ 
are all of internal edges of label $k-\delta,k,k+\delta$ in $F$.
Thus 
$G=e_1\cup e_2\cup e_1^*\cup e_2^*\cup a_{44}\cup b_{44}\cup\partial F$.

Suppose that
there exists a ring of label $k$ 
or $k+\delta$ in $F$.
This ring bounds a disk $E$ in $F$.
Now $v_1,v_2,v_3,v_4,v_5$
are all the white vertices in $F$,
and $v_1,v_2,v_3,v_4$ are
contained in the connected set 
$e_1\cup e_2\cup e_1^*\cup e_2^*\cup\partial F$ 
in
$(\Gamma_{k}\cup\Gamma_{k+\delta})\cap F$ 
containing $\partial F$.
Hence the disk $E$ contains 
at most one white vertex $v_5$.
Thus by Assumption~\ref{Ring},
the disk $E$ contains 
the white vertex $v_5$.
Hence by Lemma~\ref{ConditionRing},
we can modify the chart $\Gamma$
so that 
there is no ring of label $k$ nor 
$k+\delta$ in $F$.

Suppose that
there exists a ring of label $k-\delta$ in $F$.
This ring bounds a disk $E$ in $F$.
By
Assumption~\ref{Ring},
\begin{enumerate}
\item[(2)] 
the disk $E$ contains one of $v_3,v_4,v_5$.
\end{enumerate}
By (1),
the three white vertices $v_3,v_4,v_5$
are contained in the connected set
$e_1^*\cup e_2^*\cup a_{44}\cup b_{44}$ in 
$\Gamma_{k}\cup\Gamma_{k-\delta}$.
Hence by (2) we have that
all of $v_3,v_4,v_5$
are contained in the disk $E$
bounded by the ring of label $k-\delta$.
Thus
\begin{enumerate}
\item[(3)] $v_3,v_4,v_5\in E$.
\end{enumerate}
Let $G'$ be the connected component of 
$\Gamma_{k-2\delta}$ containing $v_5$.
By Lemma~\ref{LemmaWithTerminal}(1),
we have $w(G')\ge2$.
Thus there exists a white vertex $v_6$
in $\Gamma_{k-2\delta}$
different from $v_5$
with $v_6\in G'$.
By (1),
the vertex $v_6$ in $\Gamma_{k-2\delta}$ is
different from $v_1,v_2,v_3,v_4$.
Since $v_1,v_2,v_3,v_4,v_5$
are all the white vertices in $F$,
we have $v_6\not\in F$. 
Hence $E\subset F$ implies $v_6\not\in E$.
Since $v_5\in E$ by (3) and 
since $v_6\not\in E$,
the graph $G'$ in $\Gamma_{k-2\delta}$ 
intersects the ring $\partial E$
of label $k-\delta$.
This contradicts Condition (iv) 
of the definition of charts.
Hence there is no ring of label 
$k-\delta$ in $F$.
\end{Proof}

\section{IO-Calculation}
\label{s:IOC}

Let $\Gamma$ be a chart,
 and $v$ a vertex. 
Let $\alpha$ be a short arc of $\Gamma$ in a small neighborhood of $v$ with $v\in \partial \alpha$. 
If the arc $\alpha$ is oriented to $v$, then $\alpha$ is called {\it an inward arc}, 
and otherwise $\alpha$ is called {\it an outward arc}.

\begin{figure}
\centerline{\includegraphics{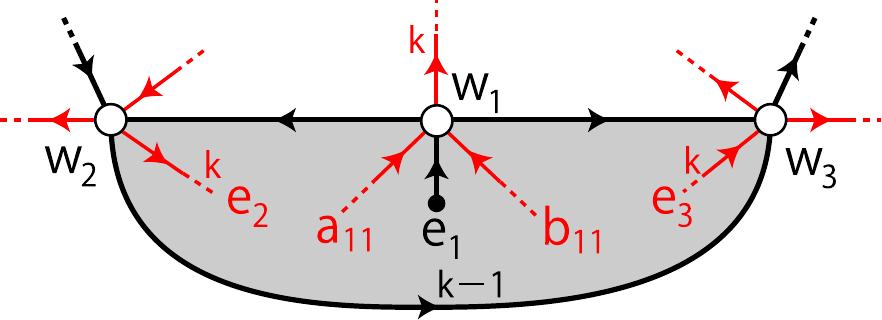}}
\caption{\label{fig18} The gray region is the disk $F$.}
\end{figure}

Let $\Gamma$ be an $n$-chart. 
Let $F$ be a closed domain with $\partial F\subset \Gamma_{k-1}\cup\Gamma_{k}\cup \Gamma_{k+1}$ for some label $k$ of $\Gamma$, where $\Gamma_0=\emptyset$ and $\Gamma_{n}=\emptyset$. 
By Condition (iii) for charts,
in a small neighborhood of each white vertex, there are three inward arcs and three outward arcs.
Also in a small neighborhood of each black vertex, there exists only one inward arc or one outward arc.
We often use the following fact, 
when we fix (inward or outward) arcs 
near white vertices and black vertices: 
\begin{enumerate}
\item[$(*)$]
{\it The number of inward arcs contained in $F\cap \Gamma_k$ is equal to the number of outward arcs in $F\cap \Gamma_k$.
}
\end{enumerate}
When we use this fact, 
we say that we use {\it IO-Calculation with respect to $\Gamma_k$ in $F$}.
For example, in a minimal chart $\Gamma$, 
consider the pseudo chart as shown in Fig.~\ref{fig18} where
\begin{enumerate}
\item[(1)] $F$ is a $3$-angled disk of $\Gamma_{k-1}$,
\item[(2)]  $w_1,w_2,w_3$ are white vertices in $\partial F$  with  $w_1,w_2,w_3\in\Gamma_{k-1}\cap\Gamma_{k}$,
\item[(3)] $e_1$ is a terminal edge of label $k-1$ containing $w_1$,
\item[(4)]  for $i=2,3$
the edge $e_i$ is of label $k$
with $w_i\in e_i\subset F$,
\item[(5)]  none of the three edges $a_{11},b_{11},e_2$ contains a middle arc at $w_1$ nor $w_2$ (by Assumption~\ref{AssumeTerminal}
none of them is a terminal edge). 
\end{enumerate}
Then we can show that $w(\Gamma_{k}\cap{\rm Int}F)\ge1$.
Suppose $w(\Gamma_{k}\cap{\rm Int}F)=0$.
If $e_3$ is a terminal edge of label $k$,
then 
by (5)
the number of inward arcs in $F\cap \Gamma_{k}$ is three,  
but the number of outward arcs in $F\cap \Gamma_{k}$ is two. 
This contradicts the fact $(*)$. 
Similarly if $e_3$ is not a terminal edge of label $k$,
then we have the same contradiction.
Thus $w(\Gamma_{k}\cap{\rm Int}F)\ge1$.
Instead of the above argument, 
we just say that 
\begin{enumerate}
\item[]
{\it we have $w(\Gamma_{k}\cap{\rm Int}F)\ge1$ 
by IO-Calculation with respect to $\Gamma_{k}$ in $F$.}
\end{enumerate}

\section{Useful Lemmata}
\label{s:Useful}

In this section,
we review useful lemmata.

\begin{lemma}$(${\rm \cite[Theorem 1.1]{ChartAppIV}}$)$
\label{LemmaNoLoop}
There is no loop in any minimal chart with exactly seven white vertices.
\end{lemma}

The following lemma will be used in 
Case (i-1) of the proof of 
Lemma~\ref{Type322Lemma5}.

\begin{figure}
\centerline{\includegraphics{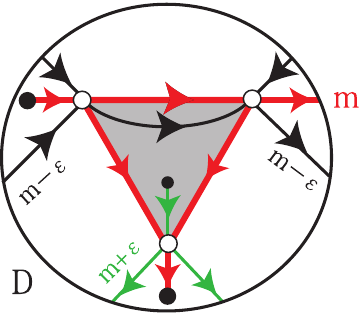}}
\caption{\label{fig19}
The gray region is the 3-angled disk $D_1$. 
The thick lines are edges of label $m$,
and $\varepsilon\in\{+1,-1\}$.}
\end{figure}

\begin{lemma}$($Triangle Lemma$)$
{\rm (\cite[Lemma 8.3(2)]{ChartAppIV})}
\label{LemmaTriangle}
 For a minimal chart $\Gamma$, 
if there exists a $3$-angled disk $D_1$ of $\Gamma_m$ without feelers in a disk $D$ as shown in Fig.~\ref{fig19},
then $w(\Gamma\cap${\rm Int}$D_1)\ge1$.
\end{lemma}

We call the graph 
in Fig.~\ref{fig04}(b)
an {\it oval}.

\begin{figure}[htb]
\centerline{\includegraphics{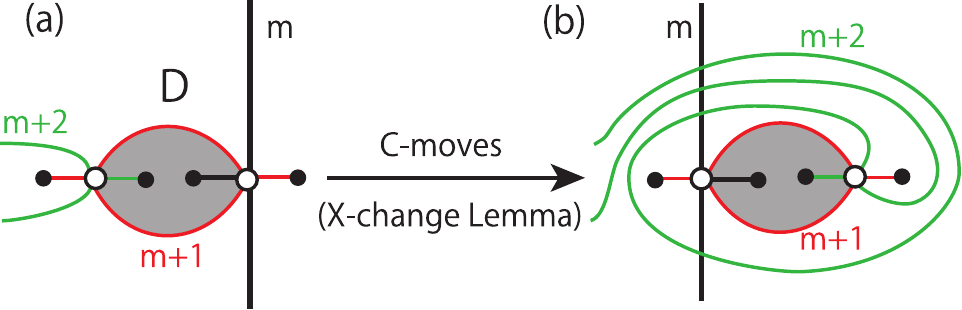}}
\caption{\label{fig20}
The gray regions are the disk $D$.}
\end{figure}

Let $\Gamma$ be a chart and $m$ a label.
An oval $G$ of $\Gamma_{m+1}$
 is said to be {\it special},
if there exists a 2-angled disk $D$ of $\Gamma_{m+1}$ without feelers
such that $\partial D\subset G$, 
$w(\Gamma\cap{\rm Int}D)=0$, 
the disk $D$ contains a terminal edge of label $m$ and a terminal edge of label $m+2$, but  $D$ does not contain any free edges, hoops nor crossings  
(see Fig.~\ref{fig20}(a)).


The following lemma will be used in 
Case (ii) of the proof of 
Lemma~\ref{Type322Lemma5}.

\begin{lemma}
{\rm (\cite[Lemma 6.1 and Lemma 6.3]{INS})}
\label{LemmaTwist}
Let $\Gamma$ be a chart.
Let $G$ be an oval of $\Gamma_{m+1}$
and $D$ a $2$-angled disk of $\Gamma_{m+1}$ without feelers such that $\partial D\subset G$ and 
$w(\Gamma\cap{\rm Int}D)=0$.
\begin{enumerate}
\item[$(1)$] $($X-change Lemma$)$
If $G$ is a special oval in 
a minimal chart $\Gamma$,
then the chart $\Gamma$ is C-move equivalent to the chart obtained from $\Gamma$ by replacing a regular neighborhood of $D$ with the pseudo chart as shown in Fig.~\ref{fig20}$($b$)$.
\item[$(2)$]
 If $D$ contains a terminal edge of label $m$ and a terminal edge of label $m+2$,
then $G$ can be modified to a special oval by C-moves in a regular neighborhood of $D$ keeping $G\cup\Gamma_m\cup\Gamma_{m+2}$ fixed.
\end{enumerate}
\end{lemma}

\section{There is no minimal chart of type $(3,2,2)$}
\label{s:MainTheorem}

In this section,
we shall show that
there is no minimal chart of type $(3,2,2)$.

The following two lemmata
will be used in the proof of 
Theorem~\ref{MainTheorem}.

\begin{figure}[htb]
\centerline{\includegraphics{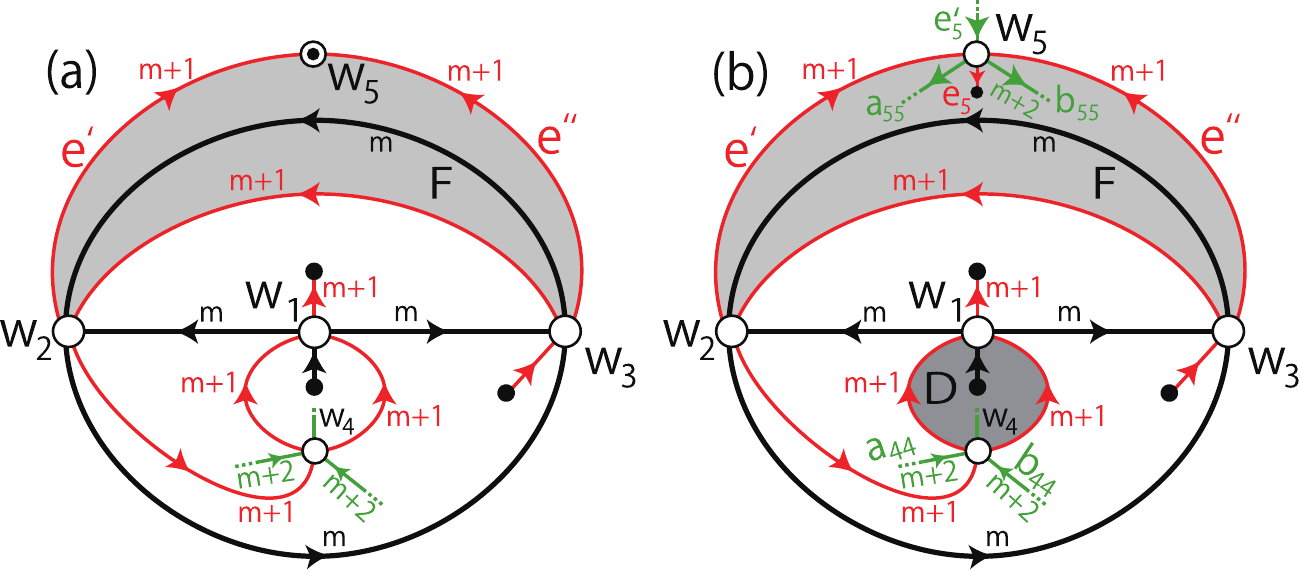}}
\caption{\label{fig21}
The gray region is the 3-angled disk $F$ of $\Gamma_{m+1}$,
and the dark gray region is the disk $D$.
}
\end{figure}

\begin{lemma}
\label{Type322Lemma5}
If a chart $\Gamma$ of type $(m;3,2,2)$
contains the pseudo chart as shown in
Fig.~\ref{fig21}$($a$)$,
then the chart $\Gamma$ is not minimal.
\end{lemma}

\begin{Proof}
Suppose that $\Gamma$ is minimal.
We use the notations as in Fig.~\ref{fig21}(a).
Here
\begin{enumerate}
\item[(1)]
$w_1,w_2,w_3\in\Gamma_m\cap\Gamma_{m+1}$,
$w_4,w_5\in\Gamma_{m+1}$.
\end{enumerate}
Since $\Gamma$ is of type $(m;3,2,2)$,
we have
\begin{enumerate}
\item[(2)]
$w(\Gamma)=7$,
\item[(3)] $w(\Gamma_m\cap\Gamma_{m+1})=3$, 
$w(\Gamma_{m+1}\cap\Gamma_{m+2})=2$,
and $w(\Gamma_{m+2}\cap\Gamma_{m+3})=2$,
\item[(4)] $w(\Gamma_{m+3})=2$.
\end{enumerate}
Thus by (1)
\begin{enumerate}
\item[(5)]
$\Gamma_{m}\cap\Gamma_{m+1}=\{w_1,w_2,w_3\}$,
$\Gamma_{m+1}\cap\Gamma_{m+2}=\{w_4,w_5\}$.
\end{enumerate}
Let $F$ be the 3-angled disk of $\Gamma_{m+1}$
with $w_2,w_3,w_5\in\partial F$
and $w_1\not\in F$.
Since
$w(\Gamma)=7$ by (2),
we have 
\begin{enumerate}
\item[(6)]
$w(\Gamma\cap${\rm Int}$F)\le2$.
\end{enumerate}
Let $e_5$ be the terminal edge at $w_5$
of label $m+1$.
Then there are two cases:
(i) $e_5\not\subset F$,
(ii) $e_5\subset F$.

{\bf Case (i).} 
Let $e_5'$ be an internal edge 
(possibly terminal edge)  
of label $m+2$ in $F$
with $w_5\in e_5'$.
By (6),
there are three cases:
(i-1) $w(\Gamma\cap${\rm Int}$F)=0$,
(i-2) $w(\Gamma\cap${\rm Int}$F)=1$,
(i-3) $w(\Gamma\cap${\rm Int}$F)=2$.

{\bf Case (i-1).}
The edge $e_5'$ must be a terminal edge. 
This contradicts Triangle Lemma 
(Lemma~\ref{LemmaTriangle}).

{\bf Case (i-2).}
There exists 
a white vertex in ${\rm Int}F$,
say $w_6$.
Thus (3) and (5) imply
$w_6\in\Gamma_{m+2}\cap \Gamma_{m+3}$.
Around the white vertex $w_6$,
there are three internal edges (possibly terminal edges)
of label $m+2$ containing $w_6$.
By the definition of charts,
only one of the three edges contains a middle arc
at $w_6$.
Hence the other two edges, say $e_6,e_6'$, 
do not contain middle arcs
at $w_6$.
By Assumption~\ref{AssumeTerminal},
neither $e_6$ nor $e_6'$ is a terminal edge.
Thus  $w(\Gamma\cap${\rm Int}$F)=1$ implies that
one of $e_6,e_6'$ contains 
the white vertex $w_5$
and the other is a loop.
This contradicts Lemma~\ref{LemmaNoLoop}.

{\bf Case (i-3).}
There exist 
two white vertices in ${\rm Int}F$,
say $w_6,w_7$.
Thus (3) and (5) imply 
$w_6,w_7\in\Gamma_{m+2}\cap \Gamma_{m+3}$.
Hence by (4)
we have 
\begin{enumerate}
\item[(7)]
$S^2-F$ does not contain any white vertices in $\Gamma_{m+3}$.
\end{enumerate}
Let $N$ be a regular neighborhood of $\partial F$ in $F$.
Let $E=Cl(F-N)$.
Then 
\begin{enumerate}
\item[(8)]
$w(\Gamma\cap{\rm Int}E)=2$ and 
$w_6,w_7\in \Gamma_{m+2}\cap\Gamma_{m+3}\cap E$, 
\item[(9)] $\Gamma_{m+2}\cap\partial E=e_5'\cap \partial E=$ one point.
\end{enumerate}
Thus applying Lemma~\ref{LemmaTwoWhiteVertices}(2) for the disk $E$,
there exist two lenses of type $(m+2,m+3)$.
This contradicts Lemma~\ref{LemmaNoLens}(2).
Hence Case (i) does not occur.

{\bf Case (ii).} 
We show $w(\Gamma_{m+2}\cap{\rm Int}F)\ge1$.
Since $e_5\subset F$
by the condition of Case (ii),
there are two internal edges $a_{55},b_{55}$
of label $m+2$ in $F$
with $w_5\in a_{55}\cap b_{55}$
(see Fig.~\ref{fig21}(b)).
Since the terminal edge $e_5$
contains a middle arc at $w_5$
by Assumption~\ref{AssumeTerminal},
neither $a_{55}$ nor $b_{55}$ 
contains a middle arc at $w_5$.
Hence by Assumption~\ref{AssumeTerminal},
\begin{enumerate}
\item[(10)]
neither $a_{55}$ nor $b_{55}$
is a terminal edge.
\end{enumerate}
Let $e',e''$ be internal edges of label $m+1$
with $\partial e'=\{w_2,w_5\}$ and
$\partial e''=\{w_3,w_5\}$.
Since the three edges $e',e'',e_5$ of 
label $m+1$ contain the white vertex $w_5$
and
since $e',e''$ are oriented inward at $w_5$,
we have that  
$e_5$ is oriented outward at $w_5$ 
(see Fig.~\ref{fig21}(b)).
Hence
\begin{enumerate}
\item[(11)]
the three consecutive edges 
$a_{55},e_5,b_{55}$
are oriented outward at $w_5$.
\end{enumerate}
Thus we have $w(\Gamma_{m+2}\cap {\rm Int}F)\ge1$
by IO-Calculation 
with respect to $\Gamma_{m+2}$ in $F$.

Let $D$ be the 2-angled disk of $\Gamma_{m+1}$
in $S^2-F$ with $w_1,w_4\in\partial D$
(see Fig.~\ref{fig21}(b)).
Next we shall show that
$w(\Gamma\cap(S^2-(F\cup D)))\ge1$.
Let $a_{44},b_{44}$
be internal edges (possibly terminal edges)
of label $m+2$ in $Cl(S^2-D)$
with $w_4\in a_{44}\cap b_{44}$.
Then 
\begin{enumerate}
\item[(12)]
 $a_{44}$ and $b_{44}$
is oriented inward at $w_4$.
\end{enumerate}
Since neither $a_{44}$ nor $b_{44}$
contain a middle arc at $w_4$, 
by Assumption~\ref{AssumeTerminal}
\begin{enumerate}
\item[(13)]
neither $a_{44}$ nor $b_{44}$
is a terminal edge.
\end{enumerate}
Let $e_5'$ be the internal edge (possibly terminal edge)
of label $m+2$ at $w_5$ with $e_5'\not\subset F$.
By (11),
\begin{enumerate}
\item[(14)]
$e_5'$ is oriented inward at $w_5$.
\end{enumerate}
Hence by (12) and (13),
we have $w(\Gamma\cap(S^2-(F\cup D)))\ge1$
by IO-Calculation with respect to $\Gamma_{m+2}$
in $Cl(S^2-(F\cup D))$.

Next we shall show that
$w(\Gamma\cap{\rm Int}F)=1$ and 
$w(\Gamma\cap(S^2-(F\cup D)))=1$.
Since $D\subset S^2-F$,
we have
$$
\begin{array}{rl}
w(\Gamma\cap (S^2-F))&
=w(\Gamma\cap (S^2-(F\cup D)))+
w(\Gamma\cap D)\\
&
\ge w(\Gamma\cap (S^2-(F\cup D)))+
w(\Gamma\cap \partial D)\\
& =
w(\Gamma\cap (S^2-(F\cup D)))+
2.
\end{array}$$
Hence by (2), 
we have 
$$
\begin{array}{rl}
7  =w(\Gamma)&=w(\Gamma\cap{\rm Int}F)+
w(\Gamma\cap\partial F)+
w(\Gamma\cap (S^2-F))\\
 & \ge w(\Gamma\cap{\rm Int}F)+3+w(\Gamma\cap (S^2-(F\cup D)))+
2.
\end{array}$$
Thus $w(\Gamma\cap{\rm Int}F)+w(\Gamma\cap (S^2-(F\cup D)))\le 2$.
Since $w(\Gamma\cap{\rm Int}F)\ge1$ and
$w(\Gamma\cap (S^2-(F\cup D)))\ge 1$,
we have
\begin{enumerate}
\item[(15)]
$w(\Gamma\cap{\rm Int}F)=1$ and
$w(\Gamma\cap (S^2-(F\cup D)))=1$.
\end{enumerate}

Next we shall show that
$e_5'$ is a terminal edge.
Suppose that
$e_5'$ is not a terminal edge.
Then $e_5'$ contains a white vertex different from $w_5$,
say $w_6$.
By (12),
we have $e_5\not\ni w_4$.
Hence $w_6\not=w_4$.
Thus (13) and $w(\Gamma\cap (S^2-(F\cup D)))=1$
imply that 
the two edges $a_{44}$ and $b_{44}$ contain
the white vertex $w_6$.
Hence by (12) and (14),
the three edges $a_{44},b_{44},e_{5}'$ of label $m+2$
are oriented outward at $w_6$.
This contradicts Condition (iii) of
the definition of charts.
Thus $e_5'$ is a terminal edge.

Finally
we shall show that
there exists an oval of label $m+2$ containing 
$a_{55},b_{55}$.
Since $w(\Gamma\cap {\rm Int}F)=1$ by (15),
the two conditions (10) and (11) imply
that the set $a_{55}\cap b_{55}$ 
contains a white vertex $w_7$
different from $w_5$
and there exists a terminal edge at $w_7$
of label $m+2$.
Hence there exists an oval $G$ of label $m+2$
 with $a_{55}\cup b_{55}\subset G$. 

Let $E$ be the 2-angled disk
of $\Gamma_{m+2}$ in $F$
with $\partial E=a_{55}\cup b_{55}$.
Since $w(\Gamma\cap{\rm Int}F)=1$
by (15),
we have 
\begin{enumerate}
\item[(16)]
$w(\Gamma\cap{\rm Int}E)=0$.
\end{enumerate}
By (3) and (5),
we have 
\begin{enumerate}
\item[(17)]
$w_7\in\Gamma_{m+2}\cap\Gamma_{m+3}$.
\end{enumerate}
By (16) and Lemma~\ref{Theorem2AngledDisk},
the disk $E$ contains one of the 
two pseudo charts as in Fig.~\ref{fig05}.
Thus by (5) and (17)
the disk $E$ contains a terminal edge 
of label $m+1$
and a terminal edge of label $m+3$.
Hence by Lemma~\ref{LemmaTwist}(2)
the oval $G$ can be modified to a special oval.
Thus by X-change Lemma 
(Lemma~\ref{LemmaTwist}(1)),
the chart $\Gamma$ is C-move equivalent to 
a minimal chart by replacing 
a regular neighborhood of $E$ with 
the pseudo chart as in Fig.~\ref{fig20}(b).
Hence the chart $\Gamma$ changes a minimal chart
satisfying the condition of Case (i)
(i.e. the terminal edge $e_5$
is not contained in $F$).
By a similar way to Case (i)
we have a contradiction.
Hence Case (ii) does not occur.

Therefore $\Gamma$ is not minimal.
\end{Proof}

\begin{figure}[htb]
\centerline{\includegraphics{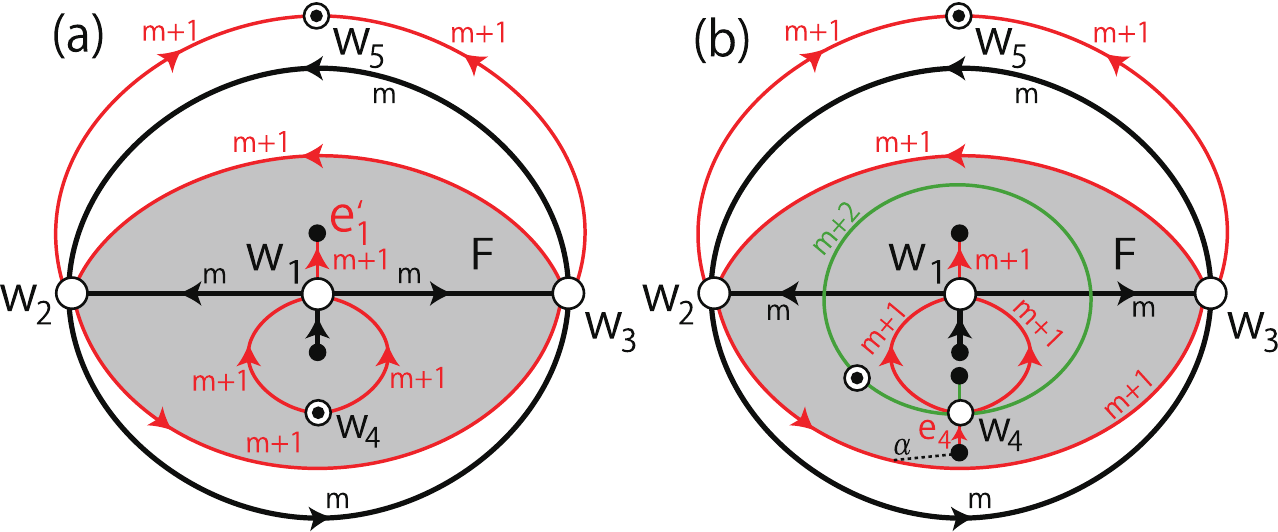}}
\caption{\label{fig22}
The gray region is the 2-angled disk $F$ of $\Gamma_{m+1}$.
}
\end{figure}

\begin{lemma}
\label{Type322Lemma6}
If a chart $\Gamma$ of type $(m;3,2,2)$
contains the pseudo chart as shown in
Fig.~\ref{fig22}$($a$)$,
then the chart $\Gamma$ is not minimal.
\end{lemma}

\begin{Proof}
Suppose that $\Gamma$ is minimal.
We use the notations as in Fig.~\ref{fig22}(a).
Here
\begin{enumerate}
\item[(1)] 
$w_1,w_2,w_3\in \Gamma_m\cap\Gamma_{m+1}$,
$w_4,w_5\in \Gamma_{m+1}$. 
\end{enumerate}
Since $\Gamma$ is of type $(m;3,2,2)$,
we have
\begin{enumerate}
\item[(2)] $w(\Gamma)=7$,
\item[(3)]
$w(\Gamma_{m}\cap\Gamma_{m+1})=3$, 
$w(\Gamma_{m+1}\cap\Gamma_{m+2})=2$, 
$w(\Gamma_{m+2}\cap\Gamma_{m+3})=2$.
\end{enumerate}
Thus by $(1)$, we have
\begin{enumerate}
\item[(4)] 
$\Gamma_{m}\cap\Gamma_{m+1}=\{w_1,w_2,w_3\}$, 
$\Gamma_{m+1}\cap\Gamma_{m+2}=\{w_4,w_5\}$. 
\end{enumerate}
Let $F$ be the $2$-angled disk of $\Gamma_{m+1}$
with $w_2,w_3\in\partial F$ and
$w_1\in {\rm Int}F$.
Then 
\begin{enumerate}
\item[(5)] 
$w_1,w_4\in {\rm Int} F$ and 
 $w_5\in S^2-F$. 
\end{enumerate}
Since $w_2,w_3\in\Gamma_m$ by $(1)$,
we have 
\begin{enumerate}
\item[(6)] 
$\Gamma_{m+2}\cap\partial F=\emptyset$.
\end{enumerate}

By $(4)$,
we have $w_4\in\Gamma_{m+2}$.
Let $G$ be a connected component of 
$\Gamma_{m+2}$
containing $w_4$.
Then by Lemma~\ref{LemmaWithTerminal}(1),
the graph $G$ contains a white vertex $w_6$ of $\Gamma_{m+2}$ different from $w_4$.
Hence by (5),
we have $G\cap{\rm Int}F\ni w_4$.
Thus (6) implies $G\subset{\rm Int}F$.
Hence $w_6\in {\rm Int}F$.
Thus by $(5)$,
we have $w_6\not=w_5$.
Since $w_6\not=w_4$,
we have 
$w_6\not\in\Gamma_{m+1}\cap\Gamma_{m+2}$
by $(4)$.
Hence $w_6\in\Gamma_{m+2}\cap\Gamma_{m+3}$.
Similarly
we can show that the open disk $S^2-F$
contains a white vertex $w_7$ of 
$\Gamma_{m+2}$ 
different from $w_5$
with $w_7\in\Gamma_{m+2}\cap\Gamma_{m+3}$.
Thus we have
\begin{enumerate}
\item[(7)] $w_6\in{\rm Int}F$
and $w_7\in S^2-F$.
\end{enumerate}
Now by (3) and (4),
we have
\begin{enumerate}
\item[(8)] $\Gamma_{m+2}\cap\Gamma_{m+3}=\{w_6,w_7\}$.
\end{enumerate}
Since $w_2,w_3\in\partial F$ and
since $w_5,w_7\in S^2-F$ by (5) and (7),
we have that
none of $w_2,w_3,w_5,w_7$ are in ${\rm Int}F$.
Hence by $(2)$,
we have $w(\Gamma\cap{\rm Int}F)\le3$.
Since $w_1,w_4,w_6\in {\rm Int}F$ by (5) and (7), 
we have 
\begin{enumerate}
\item[(9)]
$w(\Gamma\cap{\rm Int}F)=3$.
\end{enumerate}
Considering as $k=m+1,\delta=-1$,
a regular neighborhood of 
the disk $F$ contains the pseudo chart
as in Fig.~\ref{fig14}(a).
Since $w_6\in\Gamma_{m+2}\cap\Gamma_{m+3}$ by (8), 
Lemma~\ref{Lemma3PartType322}
and (9) imply that we can assume that
a regular neighborhood of $F$ 
contains one of the three pseudo charts
as in Fig.~\ref{fig15}(a),(b), (d).
Moreover 
by Lemma~\ref{Lemma2PartType322},
the chart $\Gamma$ can be modified by C-moves
keeping $G'$ fixed
so that 
\begin{enumerate}
\item[(10)]
there is no ring of 
label $m,m+1,m+2$ in $F$
\end{enumerate}
where $G'$ is the union of internal edges of
label $m,m+1,m+2$ in $F$.

{\bf Case (i).}
Suppose that a regular neighborhood of $F$ 
contains the pseudo chart
Fig.~\ref{fig15}(a) (see Fig.~\ref{fig22}(b)).
Let $e_4$ be the terminal edge at $w_4$ 
of label $m+1$.
Let $\alpha$ be an arc in $F$ connecting 
the black vertex in $e_4$
and a point in $\partial F$
with 
$G'\cap {\rm Int}\alpha=\emptyset$.
Since there is no ring of 
label $m,m+1,m+2$ in $F$ by (10), 
we can assume
$(\Gamma_m\cup\Gamma_{m+1}\cup\Gamma_{m+2})
\cap{\rm Int}\alpha=\emptyset$.
Apply C-II moves along the arc $\alpha$,
we move the black vertex in $e_4$ near 
$\partial F$. 
And we apply a C-I-M2 move between
$e_4$ and $\partial F$,
then we obtain a new terminal edge at $w_3$
of label $m+1$.
Hence $\Gamma$ contains the pseudo chart
as in Fig.~\ref{fig21}(a).
By Lemma~\ref{Type322Lemma5},
the chart $\Gamma$ is not minimal.
This contradicts 
the fact that $\Gamma$ is minimal.
Hence Case (i) does not occur.

{\bf Case (ii).}
Suppose that a regular neighborhood of $F$ 
contains the pseudo chart
Fig.~\ref{fig15}(b).
Let $e_1'$ be the terminal edge at $w_1$ 
of label $m+1$ (see Fig.~\ref{fig22}(a)).
Similarly
we can apply C-II moves so that
the black vertex in $e_1'$ is moved 
near $\partial F$.
We apply 
a C-I-M2 move between $e_1'$ and $\partial F$,
and then we obtain  
a new terminal edge of label $m+1$
containing $w_3$.
The terminal edge does not contain 
a middle arc at $w_3$.
This contradicts 
Assumption~\ref{AssumeTerminal}.
Hence Case (ii) does not occur.

{\bf Case (iii).}
Suppose that a regular neighborhood of $F$ 
contains the pseudo chart
Fig.~\ref{fig15}(d).
By a similar way to Case (ii),
we obtain a new terminal edge of label $m+1$
containing $w_3$,
and we have a contradiction.
Hence Case (iii) does not occur.

Therefore $\Gamma$ is not minimal.
\end{Proof}


\noindent{\it Proof of Theorem~\ref{MainTheorem}.}
Suppose there exists a minimal chart $\Gamma$ of type $(m;3,2,2)$.
Then $w(\Gamma_{m-1})=0$, 
$w(\Gamma_m\cap\Gamma_{m+1})=3$, and
$w(\Gamma_{m+1}\cap\Gamma_{m+2})=2$,
Thus 
\begin{enumerate}
\item[(1)] $w(\Gamma_m)=3$, 
$w(\Gamma_{m+1})=5$.
\end{enumerate}
By Lemma~\ref{LemmaNoLoop},
the chart $\Gamma$ does not contain any loop.
By Lemma~\ref{LemmaWithTerminal}(2),
the set $\Gamma_m$ contains the graph as 
in Fig.~\ref{fig04}(c).

Let $e_1$ be the terminal edge of label $m$,
and $w_1$ the white vertex in $e_1$.
We can assume that
$e_1$ is oriented inward at $w_1$
(If $e_1$ is oriented outward at $w_1$,
then 
we obtain a contradiction similarly).
Let $G$ be the connected component of $\Gamma_m$ with $w(G)=3$.
Then $G$ divides $S^2$ into three disks.
Two of them are 3-angled disks.
Let $D_1$ be the 3-angled disk
with $e_1\subset D_1$,
 $D_2$ the other 3-angled disk,
and $D_3$ the last disk
(see Fig.~\ref{fig23}(a)).
Let $w_2,w_3$ be the white vertices in $\partial D_3$
such that the internal edge $D_1\cap D_3$ is oriented from $w_2$ to $w_3$.
Considering orientation of edges around $w_3$,
the internal edge $D_2\cap D_3$ is oriented from $w_3$ to $w_2$.
Let $e_2,e_3$ be internal edges 
(possibly terminal edges) of label $m+1$ 
with $w_2\in e_2\subset D_1$ and
$w_3\in e_3\subset D_2$,
and $e',e''$ internal edges 
(possibly terminal edges) of label $m+1$
in $D_3$
containing $w_2,w_3$ respectively.
Let $a_{11},b_{11}$ be internal edges of label $m+1$ oriented inward at $w_1$
such that $a_{11},e_1,b_{11}$ lie anticlockwise around $w_1$ in this order.
If necessary we reflect the chart $\Gamma$,
we can assume that
the chart $\Gamma$ contains the pseudo chart
as in Fig.~\ref{fig23}(a).

Considering as $F=D_1$ and $k=m+1$
in the example of IO-Calculation
in Section~\ref{s:IOC},
we have $w(\Gamma_{m+1}\cap{\rm Int}D_1)\ge1$
by IO-Calculation with respect to $\Gamma_{m+1}$ in $D_1$.
Since neither $e'$ nor $e''$ contains a middle arc $w_2$ nor $w_3$, 
neither $e'$ nor $e''$ is a terminal edge
by Assumption~\ref{AssumeTerminal}. 
Thus 
by IO-Calculation with respect to $\Gamma_{m+1}$ in $D_3$,
we have $w(\Gamma_{m+1}\cap{\rm Int}D_3)\ge1$.
Since $w_1,w_2,w_3\in\Gamma_{m+1}$
and since $w(\Gamma_{m+1})=5$ by $(1)$, 
we have 
\begin{enumerate}
\item[(2)]
$w(\Gamma_{m+1}\cap{\rm Int}D_1)=1$, $w(\Gamma_{m+1}\cap{\rm Int}D_2)=0$
and $w(\Gamma_{m+1}\cap{\rm Int}D_3)=1$.
\end{enumerate}
Let $w_4$ be the white vertex in 
$\Gamma_{m+1}\cap {\rm Int }D_1$,
and
 $w_5$ the white vertex in 
$\Gamma_{m+1}\cap {\rm Int }D_3$.

Next we show that
$\Gamma$ contains the pseudo chart as in
Fig.~\ref{fig23}(b).
First look at the edge $e_3$ in the disk $D_2$.
Since the edge $e_3$ does not contain a middle arc at $w_3$,
the edge $e_3$ is not a terminal edge
by Assumption~\ref{AssumeTerminal}.
Since $w(\Gamma_{m+1}\cap{\rm Int}D_2)=0$,
we have $e_3\ni w_2$
and there exists a terminal edge $e_1'$ 
at $w_1$ of label $m+1$ in $D_2$.

Second look at the disk $D_3$.
Since neither $e'$ nor $e''$
is a terminal edge
and since $w(\Gamma\cap {\rm Int}D_3)=1$,
 both of $e',e''$
contain the white vertex $w_5$.
And there exists a terminal edge at $w_5$ 
of label $m+1$
in $D_3$.

Finally look at the disk $D_1$.
We show that $a_{11}\ni w_4$.
Since $a_{11}$ is not a terminal edge,
we have $a_{11}\ni w_2$ or $a_{11}\ni w_4$.
If $a_{11}\ni w_2$ (i.e. $a_{11}=e_2$),
then there exists a lens of type $(m,m+1)$.
This contradicts Lemma~\ref{LemmaNoLens}(2).
Hence $a_{11}\ni w_4$.
We show that $b_{11}\ni w_4$.
Since $b_{11}$ is not a terminal edge,
we have $b_{11}\ni w_2$ or $b_{11}\ni w_4$.
If $b_{11}\ni w_2$ (i.e. $b_{11}=e_2$),
then there exists a loop containing $w_4$.
This contradicts Lemma~\ref{LemmaNoLoop}.
Hence $b_{11}\ni w_4$. 
Thus 
$\Gamma$ contains the pseudo chart as in
Fig.~\ref{fig23}(b).

Since $e_2$ does not contain 
a middle arc at $w_2$, 
there are two cases:
$e_2\ni w_3$, or
 $e_2\ni w_4$.

If $e_2\ni w_3$,
then there exists a terminal edge at $w_4$
of label $m+1$.
Thus $\Gamma$ contains the pseudo chart 
as in Fig.~\ref{fig22}(a).
This contradicts Lemma~\ref{Type322Lemma6}.

If $e_2\ni w_4$,
then there exists a terminal edge at $w_3$ 
of label $m+1$
in $D_1$.
Thus $\Gamma$ contains the pseudo chart 
as in Fig.~\ref{fig21}(a).
This contradicts Lemma~\ref{Type322Lemma5}.

Therefore there is no minimal chart of type
$(3,2,2)$.
{\hfill {$\square$}\vspace{1.5em}

\begin{figure}[htb]
\centerline{\includegraphics{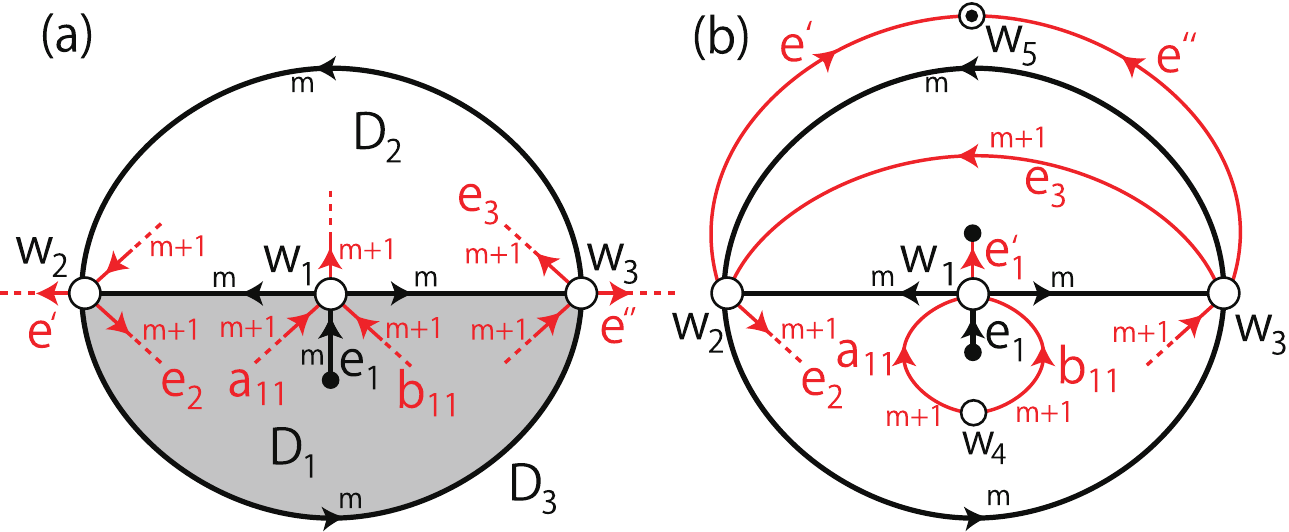}}
\caption{\label{fig23}
Charts of type $(m;3,2,2)$,
the gray region is the disk $D_1$.
}
\end{figure}




\vspace{5mm}

\begin{minipage}{65mm}
{Teruo NAGASE
\\
{\small Tokai University \\
4-1-1 Kitakaname, Hiratuka \\
Kanagawa, 259-1292 Japan\\
\\
nagase@keyaki.cc.u-tokai.ac.jp
}}
\end{minipage}
\begin{minipage}{65mm}
{Akiko SHIMA 
\\
{\small Department of Mathematics, 
\\
Tokai University
\\
4-1-1 Kitakaname, Hiratuka \\
Kanagawa, 259-1292 Japan\\
shima@keyaki.cc.u-tokai.ac.jp
}}
\end{minipage}

\end{document}